\documentclass[preprint,11pt,authoryear]{elsarticle}

\usepackage[utf8]{inputenc}
\usepackage[T1]{fontenc}
\usepackage{lmodern}
\usepackage[english]{babel}

\usepackage{booktabs}
\usepackage{graphicx}
\usepackage{amsmath}
\usepackage{amssymb}
\usepackage{amsthm}
\usepackage{euscript}
\usepackage{mathtools}
\usepackage{algorithm}
\usepackage[noEnd=false]{algpseudocodex}
\usepackage{enumerate}

\usepackage[a4paper,margin=25mm]{geometry}
\usepackage[colorlinks=true,pdfpagemode=UseOutlines,bookmarksopen,pdfstartview=FitH]{hyperref}

\biboptions{round,authoryear}
\makeatletter
\let\old@citet\citet
\renewcommand{\citet}[1]{{\renewcommand{\&}{and}\old@citet{#1}}}
\makeatother

\newtheorem{prop}{Property}

\newcommand\prob{1|p\text{-batch},r_j,s_j\leq q,\text{incompat},p_f|\sum F_j}

\journal{European Journal of Operational Research}

\makeatletter
\def\ps@pprintTitle{%
  \let\@oddhead\@empty
  \let\@evenhead\@empty
  \def\@oddfoot{%
    \normalfont\footnotesize\itshape
    Preprint submitted to \@journal\ on November 11, 2026.\hfill
  }%
  \let\@evenfoot\@oddfoot
}
\makeatother

\begin{document}

\begin{frontmatter}
\title{Parallel batch scheduling for optimizing hospital sterilization operations: arc-flow formulations, column generation and Lagrangian heuristics for flowtime minimization}

\author[cs]{Alessandro Druetto}\corref{cor1}\ead{alessandro.druetto@unito.it}
\author[math]{Andrea Grosso}\ead{andrea.grosso@unito.it}

\address[cs]{Computer Science Department, University of Turin, Corso Svizzera 185, 10149 Torino, Italy}
\address[math]{Department of Mathematics ``G.~Peano'', University of Turin, Via Carlo Alberto 10, 10123 Torino, Italy}

\cortext[cor1]{Corresponding author: Alessandro Druetto -- Email: alessandro.druetto@unito.it}

\begin{abstract}
    In a recent work by Schlembach et al., the problem of optimally managing sterilization operations in hospitals has been modeled as a machine scheduling problem with parallel batch processing and minimum flowtime objective, and heuristically solved by column generation.
    We draw upon this approach, developing arc-flow formulations and related column-generation procedures, using an exact dynamic programming algorithm to solve the pricing problems.
    We further derive Lagrangian relaxations based on minimum cost flows and shortest paths, with subgradient optimization and the volume algorithm guiding batch generation.
    A makespan bound for schedules minimizing total flow time provides a finite horizon when no operational limit is known.
    Computational experiments on literature benchmarks and more challenging generated instances demonstrate the effectiveness of the proposed approaches.
    The arc-flow model solves all literature benchmark instances to proven optimality within seconds, while the Lagrangian heuristics deliver high-quality schedules with substantially reduced computational effort.
\end{abstract}

\end{frontmatter}

\section{Introduction}
\label{sec:intro}
We consider the following scheduling problem. A set of jobs $N=\{1,2,\dots,n\}$ 
with known release dates $r_1,\dots,r_n$ and sizes $s_1,\dots,s_n$ is 
given. Such jobs must be nonpreemptively 
processed by a set of identical parallel machines $M=\{1,2,\dots,m\}$.
The job set $N$ is partitioned into a set of families $N_1,N_2,\dots,N_p$, 
let $F=\{1,\dots,p\}$. 
The machines operate by parallel batch processing: the jobs must be grouped 
into batches $B$, each containing jobs from the same family with total 
size $s(B)=\sum_{j\in B}s_j\leq q$, where $q$ is a given capacity, valid on all 
machines. Each batch of jobs from a family $f$ is processed for $p_f$ time units 
on some machine; if the processing of a batch $B\subseteq N_f$ starts at time $t$, 
all jobs $j\in B$ are started at $t$ and share the same completion time 
$C_j=t+p_f$.
The problem calls for creating the batches and scheduling them in a schedule 
$\sigma$ so that the 
total flow time $F(\sigma)=\sum_j F_j(\sigma)=\sum_{j\in N}[C_j(\sigma)-r_j]$ 
is as small as possible.
In the usual three-field notation the problem can be characterized as $\prob$.

Optimizing performance measures for batching processors is a 
well-established area within the broad field of machine scheduling. 
The early survey by \citet{Potts00} has been a solid reference 
for a number of years. Later surveys are given by \citet{Mathirajan06} and \citet{Monch11}, explicitly 
targeting the problems in the semiconductor industry. 
Finally, \citet{Fowler22} also give, to the authors' knowledge, 
the most recent survey specifically dealing with parallel batch processing.

According to \citet{Fowler22} the bulk of applications for scheduling with 
parallel batching arise in the semiconductor industry, with the management 
of diffusion ovens and burn-in operations, but the sectors of aircraft, shoes 
and automobile gear manufacturing, and health-care management deliver important 
applications as well.
This paper stems indeed from a recent work in the health-care sector 
proposed by \citet{Schlembach26}.
In the practical setting described, each surgical instrument $j$ to be sterilized 
(i.e., a job) is released at the end of some surgical operation performed during the day.
A set of $m$ identical sterilizing machines is available, where the surgical 
instruments are processed in batches that must be formed and scheduled.
Capacity constraints have to be taken into account in order to form the batches 
since each instrument $j$ has some size $s_j$ and the sterilizing machine ``trail'' 
has a limited capacity $q$.
The instruments belong to different families, and each family has its own 
requirements about the temperature and the time required for sterilization; 
specifically, each batch must be composed of jobs/instruments from the same family 
$f$ that induces a given processing time $p_f$.
\citet{Schlembach26} propose a heuristic procedure that 
relies on a very large relaxation solved by column generation.

Network-flow formulations of scheduling problems over a time-expanded graph
go back at least to \citet{PicardQueyranne78}; the
time-indexed formulations found in \citet{SousaWolsey92} admit the
same reading, made explicit by \citet{vandenAkker00},
whose columns are paths in the graph of time points. More generally, a
discrete dynamic program induces a flow-based linear formulation of its
solution set \citep{MartinRardinCampbell90}. Arc-flow models, introduced by
\citet{ValerioDeCarvalho99} for bin packing, belong to
the same family; \citet{deLima22} give a recent unified
treatment of arc-flow formulations based on dynamic programming.

Arc-flow models also lead very naturally to column-generation approaches, 
\citep[see][]{Pessoa10} where very large linear programs with exponential numbers of 
variables are efficiently handled by dynamically generating variables and columns 
using problem-specific \textit{pricing} routines. 
To the authors' knowledge, \citet{Pessoa22} and \citet{Druetto23} 
first applied ideas from arc-flow models and related column-generation approaches 
to problems with parallel batching and incompatible families.

Following such trends in the literature, the contributions of this paper are 
the followings.
\begin{itemize}
    \item An arc-flow model with a pseudopolynomial number of variables 
    is developed and tested for the problem under discussion; although the model 
    shows its limits due to the pseudopolynomial number of variables, it is able 
    to solve the whole set of real-world instances, drawn from \citet{Schlembach26},
    within seconds of computation per instance on an average PC.
    \item A realated column-generation approach is developed, with an optimized 
    pricing procedure based on dynamic programming. Instead of pursuing 
    a full branch and price approach, we follow
    \citet{Schlembach26}, using the column generation  
    in a heuristic procedure. The faster
    pricing allows to sensibly reduce the computation times.
    \item The heuristic column-generation procedure is further simplified and 
    enhanced by exploiting lagrangian relaxation techniques. This allows to 
    formulate a simple and powerful heuristic procedure, that in principle can 
    completely avoid recourse to sophisticated linear programming solvers 
    --- although we did not completely pursue this modification.
    \item Extensive computational experiments are performed on all the 
    proposed approaches.
\end{itemize}

The paper is organized as follows. In Section~\ref{sec:arc flow models} we recall 
the essentials of arc-flow models and develop the relevant variants for our 
problem: a time-indexed model and a very large model suitable for column 
generation.
In Section~\ref{sec:relaxations} two lagrangian relaxations for the very large 
model are introduced. In Section~\ref{sec:colgen} we develop different flavors 
of a heuristic procedure based on column generation. 
A routine based on dynamic programming allows for a fast pricing step.
In Sections~\ref{sec:experiments} and~\ref{sec:results} we report 
computaional experience with all the developed models and procedures.

\section{Arc-flow models for the batching problem}
\label{sec:arc flow models}
In this section we draw upon the arc-flow representations for 
scheduling models, tailoring them to our specific problem. Some known 
bound $T$ on the maximum duration of an optimal schedule is assumed to 
be known in advance.
A ``compact'' graph-based model with a pseudopolynomial number 
of arcs $\mathcal O(T n)$ and a very large model with the order of 
$\mathcal O(T 2^n)$ arcs are derived.\medskip

\noindent{\bf The basic pseudopolynomial model.}\enskip
We consider a (multi)graph $G=(V,A)$ where the set of nodes is 
$V=\{0,1,\dots,T\}$, corresponding to time slots over a sufficiently large time 
horizon within which the optimal schedule is contained. 
%We defer the definition of $T$ to Section~\ref{sec:tmax}. 
The set of arcs $A$ is defined as 
  \begin{align*}
    A&= A_1\cup A_2\cup\dots\cup A_p\cup A_0,
    \intertext{with}
    A_f&= \{(t,t+p_f;k)\mid t=0,\dots,T-p_f,\,k=1,\dots,m\}
      &&\text{for $f=1,\dots,p$}\\
    A_0&=\{(t,t+1)\mid t=0,..,T-1\}.
  \end{align*}    
The arc $(t,t+p_f;k)\in A_f$ has tail $t$, head $\tau=t+p_f$ and is labeled 
by a machine identifier $k$; it represents the possible execution of a batch 
of family $f$ over the time slot from $t$ to $t+p_f$ on machine $k$. 
Each arc $(t,t+1)\in A_0$ represents a time slot during which one or more 
machines are idle.

A relatively simple MILP can be formulated for the $\prob$ problem,
on this graph.
The scheduling of the batches is represented as the routing of $m$ units of flow 
(one unit for each machine) from the source node $t=0$ to the sink node $T$.
A binary variable $x_{tfk}$ is set to $1$ iff a batch of family $f$ is started 
at time $t$ on machine $k$, corresponding to a unit flow being routed on 
arc $(t,t+p_f;k)\in A_k$.
A binary variable $y_{jt}^k$ is set to $1$ iff a job $j$ is started at time $t$ 
in some batch on machine $k$. With these variables, the problem is modeled by the following integer linear program~\eqref{tif:obj}--\eqref{tif:vars}.
\begin{align}
    \text{minimize}\ z=&\mathrlap{\sum_{j\in N}\sum_{t=0}^{T-p_f}\sum_{k\in M} (t+p_f-r_j)y_{jt}^k}\label{tif:obj}\\
    \text{subject to}\ &\mathrlap{\sum_{k\in M}\sum_{t=r_j}^{T-p_f}y_{jt}^k=1} 
                           & & \mathllap{j\in N}\label{tif:job exec}\\
        &\mathrlap{\sum_{k\in M}\sum_{j\in N_f:\mathrlap{r_j\geq t}}s_jy_{jt}^k
        \leq q x_{tfk}} & & \mathllap{t\in V,k\in M,f\in F}
        \label{tif:batch}\\
    &x_{0,1}+\sum_{k\in M}\sum_{f=1}^px_{0fk}& & =m\label{tif: flow source}\\
    &x_{t,t+1}+\sum_{k\in M}\sum_{f=1}^px_{tfk}
               &-x_{t-1,t}-\sum_{k\in M}\sum_{f=1}^px_{t-p_f,fk}&=0
               &t=1,\dots,T-1\label{tif: flow transit}\\
    & &-x_{T-1,T}-\sum_{k\in M}\sum_{f=1}^px_{T-p_f,fk}&=-m
        \label{tif: flow sink}\\
    &\mathrlap{x_{tfk}\in\{0,1\},y_{jt}^k\in\{0,1\},\;0\leq x_{t,t+1}\leq m.}
      \label{tif:vars}
\end{align}
Constraints~\eqref{tif:job exec}--\eqref{tif:batch} deal with job execution and 
batching, requiring that \eqref{tif:job exec}~each job has a starting 
time assigned, and 
\eqref{tif:batch}~at time $t$ some job in family $f$ can be started only if a 
batch for that family starts at time $t$.
Constraints~\eqref{tif: flow source}--\eqref{tif: flow sink} are flow conservation 
constraints requiring to push $m$ units of flow from $0$ to $T$, through 
``batch processing'' arcs $(t,t+p_f;k)\in A_f$ corresponding to three-indices 
variables $x_{tfk}$ or idle arcs $(t,t+1)\in A_0$ corresponding to two-indices 
variables $x_{t,t+1}$.\medskip

\noindent{\bf A very large model.}\enskip
Let $N_f(t)=\{j\in N_f\colon t\geq r_j\}$ be the set of all jobs of family $f$ that 
can be started at time $t$.
The network $G=(V,A)$ over which the scheduling of batches can be described as a 
flow is defined by $V=\{0,1,\dots,T\}$ and set of arcs 
$A=A_1\cup\dots\cup A_f\cup A_0$ with
    \begin{align*}
        A_f=&\{(t,t+p_f;B)\mid B\subseteq N_f(t), B\neq\emptyset, 
          s(B)\leq q\}&f\in F\\
        A_0=&\{(t,t+1)\mid t=0,\dots,T-1\}.
    \end{align*}
Arcs in the $A_f$ sets are ``labeled'' with a corresponding batch $B$, we call 
them batch-arcs.
We route a unit of flow on an arc $(t,t+p_f;B)\in A_f$ when the batch $B$ is 
started at time $t$ on some machine. In a feasible schedule, idle arcs $(t,t+1)$ 
will carry a flow with value in $\{0,1,\dots,m\}$ representing the amount of 
processing resources (i.e., the number of machines) that are kept idle in the time 
slot~$t$. Note that the number of arcs in $G$ now is ``very large'' 
(exponential in the number of jobs).
The analogous of the integer linear program~\eqref{tif:obj}--\eqref{tif:vars}
is the integer program~\eqref{cg:obj}--\eqref{cg:vars} below. 
In order to ease the writing we introduce the following notations;
we define the incidence vector of a batch $B$ 
as $\boldsymbol a_B=(\text{$a_{B,j}=1$ if $j\in B$, else $0$, for $j\in N$})$;
and the all-ones column vector $\mathbf 1=(1,1,\dots,1)^T\in\{0,1\}^n$.
Also, let 
\begin{align*}
    \delta_f(t)&=\{(t,B)\mid (t,t+p_f;B)\in A_f\}\\
    \delta_f^{-1}(t)&=\{(\tau,B)\mid 
            (\tau,\tau+p_f;B)\in A_f,\tau+p_f=t\}
\end{align*} 
be the set of ($t,B$ pairs identifying) outgoing and ingoing batch-arcs of $A_f$, respectively, for node $t\in V$.

Binary variables $x_{tB}^f$ are used to represent the flow on arcs $(t,t+p_f;B)$.
The cost of starting a batch $B$ for family $f$ at time $t$ is given by
\begin{equation*}
    c_{t B}^f=\sum_{j\in B}(t+p_f-r_j)=\lvert B\rvert (t+p_f)-\sum_{j\in B}r_j.
\end{equation*}
The integer linear program for the problem writes out as follows.
    \begin{align}
        \min\ z=&\mathrlap{\sum_{f\in F}\sum_{t\in V}
          \sum_{(t,B)\mathrlap{\in\delta_f(t)}} c_{t,B}^fx_{tB}^f}\label{cg:obj}\\
        \text{subject to}\ &\mathrlap{%
                            \sum_{f\in F}
                            \sum_{(t,\tau;\mathrlap{B)\in A_f}}
                        \boldsymbol a_Bx_{tB}^f=\mathbf1}\label{cg:partition}\\
        &x_{0,1}&+\sum_{f\in F}
                 \sum_{(t,B)\mathrlap{\in\delta_f(0)}} x_{tB}^f&&&=m
                 \label{cg:flow source}  \displaybreak[2] \\ 
        &x_{t,t+1}&+\sum_{f\in F}
                   \sum_{(t,B)\mathrlap{\in\delta_f(t)}} x_{tB}^f
                   &-x_{t-1,t}
                   &-\sum_{f\in F}
                     \sum_{(\tau,B)\mathrlap{\in\delta_f^{-1}(t)}} x_{\tau B}^f
                   &=0&\quad t=1,\dots,T-1\label{cg:flow transit}\\
        &        &&-x_{T-1,T}
                    &-\sum_{f\in F}
                    \sum_{(\tau,B)\mathrlap{\in\delta_f^{-1}(T)}} x_{\tau B}^f
                   &=-m\label{cg:flow sink}\\
        &\mathrlap{x_{tB}^f\leq1, x_{t,t+1}\leq m}
        &\mathrlap{\qquad f\in F,(t,\tau,B)\in A_f,t\in V.}\label{cg:bounds}\\
        &\mathrlap{x_{tB}^f\in\mathbb Z_+,\ x_{t,t+1}\in\mathbb Z_+}
        &\mathrlap{\qquad f\in F,(t,\tau,B)\in A_f,t\in V.}\label{cg:vars}
    \end{align}
Constraints~\eqref{cg:partition} require that every job is included in 
exactly one of the scheduled batches. Constraints~\eqref{cg:flow source}--%
\eqref{cg:flow sink} are analogous to the flow-conservation 
constraints~\eqref{tif: flow source}%
--\eqref{tif: flow sink} for the time-indexed formulation. 
Note that it is no longer necessary to keep track of the machine assignments 
in the schedule: since the machines are all identical, the model simply 
ensures that no more than $m$ batches are executed at any time point, by routing 
exactly $m$ units of flow. The per-machine schedules can be recovered by 
applying a path-decomposition to the resulting flow.\medskip

\noindent{\bf Bounding the time horizon.}\enskip
Any time-indexed formulation requires some estimate of the time horizon 
over which the optimal schedule has to be deployed. 
From a practical point of view, it can be the case that some reasonable 
time horizon is naturally available for the problem. For example, in the case 
of sterilization operations in hospitals, the sterilization facility 
usually has sufficient capacity for processing practically every set of jobs 
within 24 hours. 
In some other application this is not the case, and 
in order to derive an appropriate time horizon we consider the following 
scheme.
Let $\sigma$ be any feasible schedule and denote by $C_j(\sigma)$ the 
completion time of job $j$ in schedule $\sigma$.
Recalling that every job in family $f$ has a common processing time 
$p_j=p_f$, 
define its ``waiting'' time in $\sigma$ as $W_j(\sigma)=C_j(\sigma)-p_j-r_j$, 
the time elapsed since its best possible completion time. The total waiting 
time $W(\sigma)$ in a schedule is a performance measure completely equivalent to 
the total flow time $F(\sigma)$:
\begin{align*}
   W(\sigma)&=\sum_{j\in N} W_j(\sigma)=\sum_{j\in N}[C_j(\sigma)-p_j-r_j]\\
            &=\sum_{j\in N}[C_j(\sigma)-r_j]-\sum_{j\in N}p_j
            =F(\sigma)-\underbrace{\sum_{j\in N} p_j.}_{\text{const.}}
\end{align*}
We bound the time horizon by the following property.

\begin{prop}\label{th:tmax}
Given any feasible schedule $\sigma$, the value
\begin{equation}
    T_{\mathrm{UB}} = \max_{j \in N} \left(r_j + p_j\right) + W(\sigma)
    \label{eq:tmax}
\end{equation}
is an upper bound on the makespan of every schedule that minimizes the total 
flow time.
\end{prop}
\begin{proof}
Given any optimal schedule $\sigma^\ast$ we have, by optimality of such schedule, 
$W(\sigma^\ast) \leq W(\sigma)$.
Hence, for every job $j\in N$:
\begin{align*}
    C_j(\sigma^\ast) &= r_j + p_j + W_j(\sigma^\ast)
                     \leq r_j + p_j + W(\sigma^\ast)\\ 
                     &\leq r_j+p_j+W(\sigma). %&\qquad \text{for all}\ j \in N.
\end{align*}
Taking the maximum over all jobs gives
$\displaystyle
    \max_{j\in N}\{C_j(\sigma^\ast)\} \leq 
    \max_{j \in N} \left(r_j + p_j\right) + W(\sigma).
$
\end{proof}

The bound~\eqref{eq:tmax} is particularly useful when no application-specific time 
horizon is known beforehand; the more aggressive the heuristic used to generate 
$\sigma$, the tighter the time horizon.
To the best of the authors' knowledge, this bound has not previously been proposed 
for the parallel-batch scheduling problem considered here.

%\subsection{Novelty of the computed bound}
%Without a prescribed horizon, a simple conservative alternative is a ``sum-all'' bound such as $\max_{j \in N} r_j + \sum_{j \in N} p_j$, obtained by processing all jobs individually and sequentially after the last release date.
%Such an estimate ignores the opportunities for simultaneous processing provided by batching and parallel machines, and can lead to an unnecessarily large network.

%By exploiting the total flowtime of a feasible heuristic schedule, $T^H$ can provide a substantially tighter horizon, reducing both the size of the time-expanded formulations and the number of time points examined during pricing.
%Its usefulness thus extends beyond applications with a known daily completion limit.

\section{Relaxations}
\label{sec:relaxations}
The integer program~\eqref{cg:obj}--\eqref{cg:vars} can be immediately relaxed 
to a linear program with continuous variables. We further explore a couple of 
lagrangian relaxations.\medskip

\noindent{\bf Lagrangian relaxation: minimum cost flow.}\enskip
By relaxing the constraints~\eqref{cg:partition} via lagrangian multipliers 
$\boldsymbol \pi=(\pi_i\colon i\in N)$ we get the following 
relaxation~\eqref{lagr:obj}--\eqref{lagr:vars} which corresponds to a minimum 
cost flow problem where $m$ units of flow are routed from the source node $0$ to 
the sink node $T$ on arcs with reduced costs 
$\bar c_{t B}^f=c_{t B}^f-\sum_{i\in B}\pi_i$.
\begin{align}
        L(\boldsymbol\pi)=\min_{x_{tB}^f}\ z=&\mathrlap{\sum_{f\in F}\sum_{t\in V}
        \sum_{(t,B)\mathrlap{\in\delta_f(t)}} c_{t B}^fx_{tB}^f
        +\boldsymbol\pi^{\sf{T}} \left( \mathbf{1} - 
        \sum_{f\in F} \sum_{(t,\tau;\mathrlap{B)\in A_f}}\boldsymbol a_B x_{t B}^f \right)}\notag\\
        =&\mathrlap{
           \sum_{f\in F}\sum_{t\in V}\sum_{(t,B)\mathrlap{\in\delta_f(t)}}
             \bar c_{t B}^fx_{tB}^f
             +\underbrace{\sum_{j\in N}\pi_j}_{\text{const.}}}\label{lagr:obj}\\
        \text{subject to}\quad
        &x_{0,1}&+\sum_{f\in F}
                 \sum_{(t,B)\mathrlap{\in\delta_f(0)}} x_{tB}^f&&&=m
                 \label{lagr:flow source}\\
        &x_{t,t+1}&+\sum_{f\in F}
                   \sum_{(t,B)\mathrlap{\in\delta_f(t)}} x_{tB}^f
                   &-x_{t-1,t}
                   &-\sum_{f\in F}
                     \sum_{(\tau,B)\mathrlap{\in\delta_f^{-1}(t)}} x_{\tau B}^f
                   &=0&\quad t=1,\dots,T-1\label{lagr:flow transit}\\
        &        &&-x_{T-1,T}
                    &-\sum_{f\in F}
                    \sum_{(\tau,B)\mathrlap{\in\delta_f^{-1}(T)}} x_{\tau B}^f
                   &=-m\label{lagr:flow sink}\\
        &\mathrlap{x_{tB}^f\leq1, x_{t,t+1}\leq m}
        &\mathrlap{\qquad f\in F,(t,\tau;B),t\in V.}\label{lagr:bounds}\\
        &\mathrlap{x_{tB}^f\in\mathbb Z_+,\ x_{t,t+1}\in\mathbb Z_+}
        &\mathrlap{\qquad f\in F,(t,\tau;B)\in A_f,t\in V.}\label{lagr:vars}
\end{align}
The network flow structure of the relaxed problem makes the integrality
requirements~\eqref{lagr:vars} redundant.\medskip

\noindent{\bf Lagrangian relaxation: shortest path.}\enskip 
We can note that in the integer program~\eqref{cg:obj}--\eqref{cg:vars} the 
arc capacity constraints~\eqref{cg:bounds} are actually redundant: unit upper  
bounds for batch-arc variables $x_{t B}^f\leq1$ are implied by the partition 
constraints~\eqref{cg:partition}, and the upper bounds $x_{t,t+1}\leq m$ for 
the idle arcs are implied by the routing of $m$ units of flow required by 
constraints~\eqref{cg:flow source}--\eqref{cg:flow sink}.
Dropping the capacity constraints~\eqref{cg:bounds} and relaxing the partition 
constraints~\eqref{cg:partition} with lagrangian multipliers then leads to
\begin{align}
        L(\boldsymbol\pi)=\min_{x_{tB}^f}\ z=
        &\mathrlap{
           \sum_{f\in F}\sum_{t\in V}\sum_{(t,B)\mathrlap{\in\delta_f(t)}}
             \bar c_{t B}^fx_{tB}^f
             +\underbrace{\sum_{j\in N}\pi_j}_{\text{const.}}}
             \label{lagrsp:obj}\\
    \text{subject to}\ &\text{\eqref{lagr:flow source}--\eqref{lagr:flow sink} %
        and~\eqref{lagr:vars}.\notag}
\end{align}
This latter problem amounts again to route $m$ units of flow from node $0$ to 
node $T$, but this time all $m$ units are routed together along a path with 
a minimum total reduced cost.\medskip

With both relaxations, the multipliers $\boldsymbol\pi$ can be adjusted 
by means of subgradient procedures. 
With optimal multipliers, both relaxations~\eqref{lagr:obj}--\eqref{lagr:vars} 
and~\eqref{lagrsp:obj}
deliver the same lower bound provided by the continuous relaxation 
of problem~\eqref{cg:obj}--\eqref{cg:vars}.

\section{Column-generation approaches and heuristics}
\label{sec:colgen}
We now deal with the relaxations of~\eqref{cg:obj}--\eqref{cg:vars}, 
that can be handled by column generation techniques. 
In the following, we assume for the reader some basic knowledge 
about column generation: a huge linear program with billions of columns and
variables --- the so-called master problem --- can be approximated by a 
``small'' working set of columns that defines a restricted master problem. 
This working set is dynamically extended, 
leveraging some combinatorial structure in the computation of the reduced costs.
We refer to \citet{desrosiers2026branch} for an excellent text on the subject.

\subsection{Column generation with the continuous relaxation}
\label{ssec:lpcg}
By relaxing the integrality constraints~\eqref{cg:vars} the integer 
program~\eqref{cg:obj}--\eqref{cg:vars}
becomes a very large linear program with continuous variables. 
We define a restricted master problem on a smaller subgraph of $G$, say  
$G'=(V,A')$ where $A'\subset A$. Once the relaxed problem is solved on $G'$
by the classical simplex method, we search for $x_{tB}^f$ variables 
and arcs $(t,\tau,B)\notin A'$ with the smallest (most negative) reduced 
costs, add such arcs to the restricted problem and iterate up to optimality.
If no arc with negative reduced cost is found, the optimal solution for 
the restricted problem is optimal for the full master problem.
The problem of searching for arcs with minimum reduced costs --- the 
so-called pricing phase --- offers significant combinatorial structure that 
can be leveraged.

\subsubsection{Pricing}
\label{ssec:pricing}

With the pricing procedure we use dual information from the optimally solved 
restricted problem in order to detect (at least) an arc $(t,\tau;B)\in A_f$ and 
the associated variable $x_{tB}^f$ with minimum reduced cost.
This pricing problem can be
decomposed into $p$ subproblems, one per family, each asking for the computation of 
a minimum reduced cost batch $B$ for family $f$, executed at $t$. Consider column vectors 
of dual values $\boldsymbol{\pi} = (\pi_j\colon j \in N)$ and 
$\boldsymbol{\lambda} = (\lambda_t\colon t=0,\dots,T)$ for
constraints~\eqref{cg:partition} and~\eqref{cg:flow source}--\eqref{cg:flow sink}
respectively, obtained with the optimal basis for the restricted master problem.
By the well known theory of linear programming, the reduced cost for 
a variable/arc representing a batch $B$
from family $f$ to be executed at time $t$ is defined as 
\begin{equation}
  \hat{c}_{t B}^f = \sum_{j \in B} (t + p_f - r_j) - \sum_{j \in B} \pi_j -
  (\lambda_t - \lambda_{\tau}) \quad \quad \text{with $\tau = t + p_f$.}
  \label{eq:redcost}
\end{equation}
The problem of searching for the batch that minimizes the reduced
cost~\eqref{eq:redcost} is, for a given time point $t$,
\begin{align}
    \min_B\ 
      &\left\{\hat c_{t B}^f=\sum_{j\in B}(t + p_f - r_j-\pi_j)
          -(\lambda_t-\lambda_{\tau})\colon 
                           s(B)\leq q,B\subseteq N_f(t)\right\}\notag\\
     =&\min_B
       \left\{\bar\zeta(t)=\sum_{j\in B}(t + p_f - r_j-\pi_j)
          \colon s(B)\leq q, B\subseteq N_f(t)\right\}
            -(\lambda_t-\lambda_{\tau})\notag\\
    =&-\max_B\left\{\zeta(t)=\sum_{j\in B} (\pi_j+r_j-t-p_f)\colon 
          s(B)\leq q, B\subseteq N_f(t) \right\}
    -(\lambda_t-\lambda_{\tau}).
            \label{eq:pricing}
\end{align}
The maximization problem in~\eqref{eq:pricing} is a knapsack problem on the 
jobs in $N_f(t)$ that are available for execution at time $t$.
For fixed $t$ the difference $(\lambda_t - \lambda_{\tau})$ is
constant, and the contribution $-t-p_f$ to the knapsack profit is the same for 
all jobs in $B$, so solving~\eqref{eq:pricing} amounts to solving the peculiar 
knapsack problem 
\begin{equation}
          \max_B
            \left\{\zeta(t)=\sum_{j\in B}(\pi_j+r_j)-\lvert B\rvert(t+p_f)\colon
                    s(B)\leq q, B\subseteq N_f(t)\right\}
              \label{eq:knap}
\end{equation}
where part of the objective is determined by the cardinality of the selected 
batch $B$.
We show that problem~\eqref{eq:knap} can
be solved by filling a suitable dynamic programming table, and such table 
can be used for all time points $t$. 
Assume the jobs in $N_f$ are consecutively indexed by $1, 2, \ldots
n_f$ in decreasing order of release date: $r_1 \geq r_2 \geq \cdots \geq
r_{n_f}$;
this ensures that if job $i$ is available at some time $t$, so are also the 
jobs $i+1,..,n_f$. 
The value $\varphi_i (k, s)$ represents the maximum value of the total $(\pi_j
+ r_j)$ that can be obtained by selecting within a given total size $s$, $k$ jobs
from those in $N_f$ that have a release time $r_j \leq r_i$. The values of
$\varphi_i (k, s)$ for $i \in N_f$, $k = 1, \ldots, n_f$ and $s \leq q$ can be
computed by dynamic programming:
\begin{equation}
  \varphi_i (k, s) = \max \left\{\begin{array}{ll}
    \varphi_{i + 1} (k, s) & \qquad (y_i = 0)\\
    (r_j+\pi_j)+\varphi_{i + 1} (k - 1, s - v_i) & \qquad (y_i = 1)
  \end{array}\right. \label{pricing dp}
\end{equation}
with initial conditions
\begin{equation}
  \begin{array}{llll}
    \varphi_i (k, 0) & = & - \infty & \qquad \text{for all $k > 0$}\\
    \varphi_i (0, s) & = & 0 & \qquad \text{for all $s$},\\
    \varphi_{n_f}(1,s)&=&\begin{cases}
        (\pi_{n_f}+r_{n_f})&\ \text{if $s_{n_f}\leq s$}\\
        -\infty&\ \text{otherwise.}
    \end{cases}&\qquad \text{for all $s$.}
  \end{array} \label{pricing init}
\end{equation}
Given some time $t$, select the smallest index $i \in N_f$ such that $r_i \leq
t$; the optimal value for~\eqref{eq:knap} is
\[\zeta^\ast(t)=\max_k \left\{ \varphi_i (k, q) - k (t + p_f) : \; 
  \text{$k$} = 1, \ldots, n_f - i + 1 \right\}\]
and the value for the minimum reduced cost corresponding to the optimal 
value of~\eqref{eq:pricing} is given by
\begin{equation}
  \hat{c}_ {t B^\ast} = -\zeta^\ast(t)-(\lambda_t-\lambda_{\tau}).
\end{equation}
The content of the optimal batch can be recovered by backtracking on the
dynamic programming table. The pricing scheme for a family is sketched 
in Algorithm~\ref{algo:pricing1}, that returns for a given family $f$ a set 
of $(t,B)$ pairs identifying arcs and variables with negative reduced costs.
The algorithm uses the same dynamic programming table for all values of $t$, 
and returns all the negative-cost arcs it finds, not only the most negative.
If Algorithm~\ref{algo:pricing1} returns empty sets for all families, this 
certifies that the current optimal solution for the restricted 
master problem is also optimal for the whole master problem.

\begin{algorithm}[t]
    \caption{Pricing for family $f\in F$. 
        $N_f$ is sorted and indexed by nonincreasing release
        times $r_1 \geq r_2 \geq \cdots \geq r_{n_f}$.
    }\label{algo:pricing1}
    \begin{algorithmic}[1]
        \Require dual values $\boldsymbol\pi=(\pi_i:i \in N_f)$ and 
        $\boldsymbol\lambda=(\lambda_t:t=0,\dots,T)$.
        \Ensure a (possibly empty) set of arcs $S$ with negative reduced costs.
      \Function{PriceArcs}{$f$: family; $\boldsymbol\pi,\boldsymbol\lambda$: duals}
        \State Set $S:=\emptyset$
        \State Compute $\varphi_i (k, s)$ for all $k$, $s$ and $i \in N_f$.
        \For{$t = r_{n_f}, \ldots T-p_f$}
            \State Select the smallest index $i \in N_f$ such that $t \geq r_i$.
            \For{$k = 1, \ldots n_f - i + 1$}
                \State Compute the reduced arc cost 
                \Statex\qquad $
                  \hat{c} (t, B) = - \varphi_i (k, q) + k (t + p_f) 
                  - (\lambda_t -\lambda_{\tau})
                $\Comment{$\tau=t+p_f$}
                \If{$\hat{c}_{tB}<0$}
                    \State Recover $B$ by backtracking and set $S:=S\cup\{(t,\tau;B)\}$
                \EndIf
            \EndFor
        \EndFor
        \State \Return $S$
      \EndFunction
    \end{algorithmic}
\end{algorithm}

The pricing procedure for family $f$ requires at most 
$\mathcal O(q\lvert N_f\rvert^2)$ operations to fill the dynamic programming table 
(line~3) and at most $\mathcal O(T \lvert N_f\rvert)$ operations for computing 
the reduced costs. 
A complete computation of 
all $\varphi_i(k,s)$ can be avoided by starting with an empty table and 
adopting a memoization strategy, so that 
the table is filled top-down as needed, computing only the values that are 
strictly necessary.

\subsubsection{Overall heuristic procedure}
\label{ssec:heur}
A heuristic procedure for problem $\prob$ can be defined by three steps 
--- see also Algorithm~\ref{algo:overall procedure}. 
It follows the approach that is sometimes called ``price and branch'': 
solve the relaxation by column generation, then solve the final restricted 
problem as an integer one. This usually deliver very good solutions 
with a narrow certified optimality gap, 
while avoiding the intricacies of a full branch and price exact approach. 
\begin{enumerate}
    \item Compute the heuristic solution by some greedy procedure, defining 
    an initial set of scheduled batches and the corresponding arcs in the 
    restricted master problem.
    \item Iteratively solve the relaxed restricted master problem, generating 
    new arcs/batches by Algorithm~\ref{algo:pricing1} and expanding the problem, 
    until no more arcs with negative reduced costs can be found and the optimal 
    solution of the master relaxation is available.
    \item Revert to integer variables for the final restricted master problem, 
    and solve it (via exact or heuristic procedures).
\end{enumerate}
In our implementation Step~1 is performed by applying the ``batch whenever idle'' 
heuristic by \citet{Dobson01}. 
Step~3 is performed by solving the resulting integer program to optimality 
by a solver. 
The whole procedure is similar to the one designed by \citet{Schlembach26}, except for the pricing based on dynamic programming 
and the initial heuristic.
\begin{algorithm}[t]
    \caption{Overall solution procedure.}\label{algo:overall procedure}
    \begin{algorithmic}[1]
        \State{Generate a heuristic solution $\sigma$}
        \State{Set $V:=\{0,1,\dots,T\}$}
          \Comment{Initialize reduced graph $G=(V,A')$}
        \State{Set $A'$ to the set of all $(t\tau;B)$ arc-batches induced 
          by $\sigma$ plus the idle arcs in $A_0$}
        \Repeat
          \State{Solve the continuous relaxation~\eqref{cg:obj}--\eqref{cg:vars}
            restricted to $G'=(V,A')$}
          \State{Retrieve multipliers $\boldsymbol\pi$, $\boldsymbol\lambda$}
            \Comment{For constraints~\eqref{cg:partition} 
            and~\eqref{cg:flow source}--\eqref{cg:flow sink} resp.}
          \State{Set $S:=\emptyset$}
          \For{$f=1,2,\dots,p$}\Comment{Price out arcs for each family}
            \State{Set $S:=S\cup\textsc{PriceArcs}(f,\boldsymbol\pi,
            \boldsymbol\lambda$)}\Comment{By Algorithm~\ref{algo:pricing1}}
          \EndFor
          \State{Set $A':=A'\cup S$}
        \Until $S=\emptyset$\Comment{Relaxation solved at optimality}
        \State{Solve~\eqref{cg:obj}--\eqref{cg:vars} on the restricted
          arc set $A'$}\Comment{Solve MILP, getting a heuristic solution}
    \end{algorithmic}
\end{algorithm}

\subsection{Column generation with the minimum cost flow relaxation}
\label{ssec:cgmcf}
As an alternative to the continuous relaxation of the integer program%
~\eqref{cg:obj}--\eqref{cg:vars} we can consider the minimum cost flow 
problem~\eqref{lagr:obj}--\eqref{lagr:vars} arising from the lagrangian 
relaxation of the partition constraints~\eqref{cg:partition}. 
If we consider the restricted master problem on a working set $A'\subset A$, 
for any given multiplier vector $\boldsymbol\pi$, the optimal value 
$L(\boldsymbol\pi)$ can be efficiently computed using a network simplex 
routine. The multipliers $\boldsymbol\pi$ are adjusted by a subgradient 
procedure in order to maximize the lower bound $L(\boldsymbol\pi)$.
The optimal $\boldsymbol\pi$ and the corresponding 
flow potentials $\boldsymbol\lambda^\ast$ will give the 
optimal dual variables for constraints~\eqref{cg:partition}--\eqref{cg:flow sink}.
Hence, solving the lagrangian dual of a restricted master problem 
is a viable alternative to invoking the simplex method.

An iterative subgradient optimization procedure can be applied in order to 
determine the optimal multipliers.
We deliberately adopted a very basic implementation, for two reasons. First, 
given the network flow nature of the lagrangian subproblem, each
iteration is cheap. Second, we do not need the master relaxation solved to high
accuracy, only dual information good enough to steer the pricing; in this regime
properly tuned subgradient methods are considered competitive with more
sophisticated approaches \citep{Frangioni17}. 
%Third, the Lagrangian value is a valid lower bound
%at every iteration \citep{Geoffrion74}, unlike the value of the restricted master.
Any more refined schemes can replace the basic one without changing the
rest of the procedure.

Given the current solution for~\eqref{lagr:obj}--\eqref{lagr:vars}, a subgradient 
vector $\boldsymbol g$ is computed and the multipliers $\boldsymbol\pi$ 
after the $r$-th iteration are updated to
\[\boldsymbol\pi^{+}=\boldsymbol\pi+\alpha_r \boldsymbol g,\quad
\text{with}\quad
\boldsymbol g=\mathbf1
  -\sum_{f\in F}\sum_{(t,\tau,\mathrlap{B))\in A_f}}\boldsymbol a_Bx_{tB}.\]
Initial multiplier values are given by the flowtime values obtained by the Dobson-Nambimadom heuristic of Step 1.
The perturbation step $\alpha_r$ is set at $sz/\sqrt{r}$, with initial step size $sz = 1$.
The subgradient procedure is stopped after $\lVert\boldsymbol g\rVert$ falls 
below a tolerance of $1$ (because of the integrality of all $\boldsymbol g$ components, this means $\lVert\boldsymbol g\rVert = 0$ in practice), or after 20 consecutive iterations fail 
to increase the optimal value of~\eqref{lagr:obj}. 
As a failsafe, no more than 1000 total subgradient iterations are allowed.
Since the subgradient update results only in a perturbation of the costs 
$\bar c_{t B}^f$, the network simplex algorithm is a reasonable choice for 
handling the incremental reoptimizations.
At convergence the multipliers $\boldsymbol\pi$ will converge to the values 
of the optimal dual variables for the relaxed constraints~\eqref{cg:partition}, 
and the components $\pi_1,\dots,\pi_n$
together with the optimal minimum cost flow potentials $\lambda_t$ ($t\in V$) 
can be fed to Algorithm~\ref{algo:pricing1} to price out more arcs.
Note that using an approximate dual solution substantially voids the claim of 
optimality that can be done on the relaxed master problem in 
Algorithm~\ref{algo:overall procedure}, but this does not invalidate the 
heuristic approach.
In such a case, a proper Lagrangian bound can be computed using the relaxation described in~\eqref{lagrsp:obj}; see next section for details.

\subsection{Column generation with the shortest-path relaxation}
\label{ssec:lagrsp}
Solving relaxation~\eqref{lagrsp:obj} requires computing a shortest 
path from node $0$ to node $T$ on the whole graph $G=(V,A)$. 
The graph $G$ is acyclic, with $0,1,\dots,T$ being a valid topological ordering  
for the nodes in $V$.
Hence, the shortest path can be computed by solving a set of Bellman-like 
equations
\begin{equation}\label{eq:Bellman}
    \ell_t=\min\begin{cases}
        \ell_{t-1}\\
        \min\{\ell_{t'}+\bar c_{t'B}^f\colon (t'B)\in\delta_f^{-1}(t),f\in F\}
    \end{cases}\qquad t=0,1,\dots,T,
\end{equation}
where $\ell_t$ represents the cost of a shortest path from $0$ to $t$.
Algorithm~\ref{algo:Gondran-Minoux} solves the above equations following 
the topological ordering, sequentially updating for 
each $t=0,1,\dots,T$ and $f\in F$,
\begin{align}
    \ell_{\tau}&:=\min\{\ell_{\tau},\ell_t\}&&\quad\text{with $\tau=t+1$,}
    \label{eq:idle arc}\\
    \ell_{\tau}&:=\min\big\{
      \ell_{\tau},
      \ell_t+\min\{\bar c_{t B}^f\colon s(B)\leq q,  
                          B\subseteq N_f(t)\}
      &&\quad\text{with $\tau=t+p_f$.}\label{eq:batch arc}
\end{align}
Update~\eqref{eq:idle arc} takes care of idle arcs $(t,\tau=t+1)\in A_0$, while 
update~\eqref{eq:batch arc} takes care of batch processing arcs 
$(t,\tau=t+p_f;B)\in A_f$. This is substantially an implementation of the 
well-known Gondran-Minoux algorithm for shortest paths on acyclic graphs.
In~\eqref{eq:batch arc}, the inner minimization can 
be handled by dynamic programming, since determining a batch arc-for family $f$ 
with minimum reduced cost requires solving
\begin{align}
    &\min_B\left\{\bar c_{t B}^f=\sum_{j\in B}(t+p_f-r_j-\pi_j)\colon
             \sum_{i\in B}s_i\leq q, B\subseteq N_f(t)\right\}\notag \\
    &=-\max_B\left\{\zeta(t)=\sum_{j\in B}(\pi_j+r_j-t-p_f)\colon
             \sum_{i\in B}s_i\leq q, B\subseteq N_f(t)\right\}\notag \\
     &=-\max_B
            \left\{\zeta(t)=\sum_{j\in B}(\pi_j+r_j)-\lvert B\rvert(t+p_f)\colon
                    s(B)\leq q, B\subseteq N_f(t)\right\}.
                    \label{eq:knap1}
\end{align}             
This is the knapsack problem~\eqref{eq:knap} solved in Section~\ref{ssec:pricing}. 

\begin{algorithm}[htb]
    \caption{Solution of the shortest path 
      relaxation~\eqref{lagrsp:obj}.}
      \label{algo:Gondran-Minoux}
    \begin{algorithmic}[1]
        \Require{Multipliers $\boldsymbol\pi=(\pi_j\colon j\in N)$ 
          and job/families data.}
        \Ensure{A shortest path from $0$ to $T$, along with the set of generated 
        minimum cost arcs.}
        \State{Compute the tables $\varphi_i(k,s)$ for each family $f\in F$}
        \State{Set $\ell_0:=0$, $\ell_t=\infty$ for $t\in V\setminus\{0\}$}
        \For{$t=0,\dots,T-1$}
          \If{$\ell_t<\ell_{t+1}$}
            \State{Set $\ell_{t+1}:=\ell_t$}
              \Comment{Handles the idle-arc case}
          \EndIf
          \For{$f\in F$}
            \State{Consider $\tau:=t+p_f$}
            \If{$\tau\leq T$}
              \State{Select the smallest index $i\in N_f$ such that $t\geq r_i$}
              \For{$k=1,\dots,n_f-i+1$}%
                \Comment{Get cheapest $B$ for each $|B|=k$}
                \State{Using the $\varphi$ table for $f$ compute
                  $\bar c_{t B}^f=-\varphi_i(k,q)+k(t+p_f)$
                }
                \If{$B\neq\emptyset$ and $\ell_t+\bar c_{t B}^f<\ell_\tau$}
                  \State{Set $\ell_\tau:=\ell_t+\bar c_{t B}^f$}
                \EndIf
              \EndFor
            \EndIf
          \EndFor
        \EndFor
        \State{\Return 
          $L(\boldsymbol\pi)=\sum_j\pi_j+m \ell_T$}
    \end{algorithmic}
\end{algorithm}

The lagrange multipliers in the relaxation are iteratively adjusted by the volume
algorithm by \citet{Barahona00}, 
run for a maximum of 1000 iterations. A simpler subgradient procedure was tried
like the one illustrated in Section~\ref{ssec:cgmcf}, but it exhibited a very 
low convergence speed.

The adopted implementation of the volume algorithm is sketched in 
Algorithm~\ref{algo:volume}, with a common choice of parameters.
We refer to \citet{Barahona00} and related literature for a detailed discussion 
of the algorithm.
The step size $\Delta$ is determined on line~9 with a Polyak-style 
formula, where the value of the upper bound $\text{UB}$ is given by the 
Dobson-Nambimadom heuristic.

\begin{algorithm}[htb]
    \caption{Volume algorithm --- multiplier adjustment for the 
    shortest-path relaxation~\eqref{lagrsp:obj}. Tolerances are set
    at $\varepsilon_1 = 10^{-2}$ and $\varepsilon_2 = 10^{-4}$.}
    \label{algo:volume}
    \begin{algorithmic}[1]
        \State{Set $\boldsymbol{\bar\pi}:=\boldsymbol\pi$}
        \State{Solve~\eqref{lagrsp:obj} get $L(\boldsymbol{\bar \pi})$ and the corresponding  
          $\boldsymbol x(\boldsymbol{\bar\pi})$, 
          set $\boldsymbol {\bar x}:=\boldsymbol x(\boldsymbol{\bar\pi})$}
        \State{Set $\rho:=0.1$, $\alpha_{\max}:=0.1$}
        \For{$\texttt{iter}=1,2,\dots,1000$}
          \State{Compute subgradient $\displaystyle\boldsymbol{\bar g}=\mathbf1
  -\sum_{f\in F}\sum_{(t,\tau,\mathrlap{B))\in A_f}}\boldsymbol a_B\bar x_{tB}$
    and $\displaystyle\bar z=\sum_{f\in F}\sum_{t\in V}
            \sum_{(t,B)\mathrlap{\in\delta_f(t)}} c_{t B}^f\bar x_{tB}^f$
            }
          \If{$\lVert\boldsymbol{\bar g}\rVert_{\infty}\leq\varepsilon_1$ and 
              $\lvert \bar z-L(\boldsymbol{\bar\pi})\rvert\leq 
                \varepsilon_2\lvert L(\boldsymbol{\bar\pi}\rvert)$}
                    \State{\Return $\boldsymbol{\bar\pi}$ }
          \EndIf
          \State{Compute step size 
            $\Delta=\rho [\text{UB}-L(\boldsymbol{\bar\pi})]/
                                \lVert\boldsymbol{\bar g}\rVert^2$}
          \State{Set $\boldsymbol\pi:=\boldsymbol{\bar\pi}
                  +\Delta\cdot\boldsymbol{\bar g}$}
          \State{Solve~\eqref{lagrsp:obj}, get $L(\boldsymbol{\pi})$ and the corresponding $\boldsymbol x(\boldsymbol\pi)$}
          \State{Compute subgradient $\displaystyle\boldsymbol g=\mathbf1
  -\sum_{f\in F}\sum_{(t,\tau,\mathrlap{B))\in A_f}}\boldsymbol a_Bx_{tB}$}
          \State{Select  
          $\displaystyle\alpha^\ast:=\arg\min\{\lVert\alpha\boldsymbol g
                             +(1-\alpha)\boldsymbol{\bar g}\rVert^2\colon
                             \alpha\in[0.1\alpha_{\max},\alpha_{\max}]$\}}
            \State{Set 
                $\boldsymbol{\bar x}:=\alpha^\ast \boldsymbol x(\boldsymbol\pi)+
                                      (1-\alpha^\ast)\boldsymbol{\bar x}$}
        \If{$L(\boldsymbol\pi)>L(\boldsymbol{\bar\pi})$}
          \State{Expand $\rho:=\min\{2,1.1\rho\}$ 
            if $\boldsymbol g^{\sf T}\boldsymbol{\bar g}\geq0$}
          \State{Set $\boldsymbol{\bar\pi}:=\boldsymbol\pi$}
        \Else
          \State{Reduce $\rho:=\max\{10^{-4},0.66\rho\}$ 
            if $20$ consecutive iterations failed to increase $L$}
        \EndIf
        \State{Reduce $\alpha_{\max}:=\max\{10^{-5},\frac12\alpha_{\max}\}$ 
          if $L(\boldsymbol{\bar\pi})$ has increased less than $0.1\%$ for 
          $100$ iterations}
        \EndFor
        \State{\Return $\boldsymbol{\bar\pi}, L(\boldsymbol{\bar\pi})$}
    \end{algorithmic}
\end{algorithm}

The overall heuristic procedure of Section~\ref{ssec:heur} is now slightly 
modified as follows.
\begin{enumerate}
    \item Compute the heuristic solution by the heuristic procedure of 
    Dobson and Nambimadom.
    \item Apply the volume algorithm and Algorithm~\ref{algo:Gondran-Minoux} 
    to solve relaxation~\eqref{lagrsp:obj} and determine the best multipliers 
    (a restricted master problem is no longer needed). 
    \item Solve the integer program~\eqref{cg:obj}--\eqref{cg:vars}, restricting 
    the batch-arcs in graph $G$ to the collection of those generated in steps~1
    and~2.
\end{enumerate}

The volume algorithm is well known for being able to deliver almost-feasible 
primal solutions, with a very low violation of the relaxed primal constraints.
This seems to be the case for the final average solution 
$\boldsymbol{\bar x}$ that Algorithm~\ref{algo:volume} computes as 
a by-product. Experimentally, this $\boldsymbol{\bar x}$ is characterized by a 
very small violation of the partition constraints. 
This probably paves the way for some rounding heuristic 
that could replace the integer optimization in the last step of the solution 
procedure, but for now we did not pursue this possibility.

We consider the procedure built around relaxation~\eqref{lagrsp:obj} as the 
most interesting one.
Experimentally, it proved to be fast and effective; it is simple, since 
it does not necessarily rely on linear programming libraries.
The whole program for solving the relaxation, in its most basic form, can fit into 
a couple of hundreds lines of code.

The arc generation via dynamic programming is analogous to the pricing 
procedure described in Section~\ref{ssec:pricing}, with the relevant difference 
that here no exogenous node duals $\lambda_t$ have to be taken into account 
and the knapsack problem alone determines the arc to be generated.
This offers some leverage to optimize the generation procedure.

\begin{prop}\label{prop:relevance}
  In order to solve the shortest path equations~\eqref{eq:Bellman}, 
  the per-family optimization~\eqref{eq:knap1} 
  for a given $t$ can be restricted to batches $B\subseteq\hat N_f(t)$ where
  \[\hat N_f(t)=\{r_j\leq t<r_j+\pi_j-p_f\}\subseteq N_f(t).\]
\end{prop}
We define the interval $[r_j,r_j+\pi_j-p_f)$ as the relevance interval 
of job $j$.
\begin{proof}
    The requirement $r_j\leq t$ in the definition of $\hat N_f(t)$ is obvious. 
    We deal with the second part of the defining inequality.
    \begin{itemize}
        \item The optimization problem~\eqref{eq:knap1} aims to perform the 
        update~\eqref{eq:batch arc} which in turn is instrumental to 
        solve the shortest-path equations~\eqref{eq:Bellman}.
        \item Due to the structure of the graph, between any two time-nodes 
        $t$ and $\tau=t+p_f$ a zero-cost path from $t$ to $\tau$ is always 
        available through idle arcs $(t,t+1),(t+1,t+2),\dots,(\tau-1,\tau)$.
        \item Hence, a batch-arc with $(t,B)\in\delta_f(t)$ connecting $t$ and $
        \tau$ is relevant for the shortest path only if $\bar c_{t B}^f<0$, 
        otherwise that batch-arc can be safely ignored.
    \end{itemize}
    Thus, only solutions of problem~\eqref{eq:knap1} with $\bar c_{t B}^f<0$ 
    need to be considered. Furthermore, for problem~\eqref{eq:knap1} there 
    exists an optimal solution $B^\ast$ 
    such that no job $j\in B^\ast$ has $t+p_f-r_j-\pi_j\geq0$: otherwise, removing 
    such $j$ from $B^\ast$ would reduce $s(B^\ast)$ without increasing 
    $\bar c_{t B^\ast}^f$. Hence, we can form the batches only with jobs $j$
    such that $t+p_f-r_j-\pi_j<0\iff t<r_j+\pi_j-p_f.$
\end{proof}

Thanks to Property~\ref{prop:relevance}, for a given $t$ and family $f$ 
only the relevant jobs in $\hat N_f(t)$ need to be taken into account to solve 
problem~\eqref{eq:knap1}. This can sensibly restrict the sizes of the dynamic 
programming tables $\varphi$ used for each family.
We can consider the time points $t_1\leq t_2\leq\dots\leq t_K$ obtained by 
sorting the set $\{r_j\colon j\in N_f\}\cup\{r_j+\pi_j-p_f:j\in N_f\}$.
The time horizon $0,\dots,T$ can be split into segments 
$[0,t_1),[t_1,t_2),\dots,[t_K,T]$; within each segment, the set of 
relevant jobs for family $f$ $\hat N_f(t)$ remains constant. Hence, instead 
of keeping one copy of a $\varphi$ table for all $t$ we reinitialize it  
at each breakpoint. This offers a significant computational advantage  
when the sets $\hat N_f(t)$ tend to be small.

%\[\alpha_k=\theta_k(\mathrm{UB}-
%L(\boldsymbol\pi_k))/\lVert\boldsymbol g_k\rVert^2\]

\section{Experimental setting}
\label{sec:experiments}
We consider three sets of instances: the experimental instances and the full real-world instances of \cite{Schlembach26}, denoted by \texttt{orig} and \texttt{real}, respectively, and an additional generated set, denoted by \texttt{gen}, designed to stress batch formation and machine
contention.
Table~\ref{tab:instance-sets} summarizes their main parameters.

\begin{table}[htbp]
\centering
\begin{tabular}{l|r|lrrr}
    \toprule
    Set & Instances & $n$ & $m$ & $p$ & $q$\\
    \midrule
    \texttt{orig} & 150 & $\{10, 25, 50, 75, 100\}$ & $2, 3, 4$ & 2 & 180\\
    \texttt{real} & 10  & 181--214                  & 4         & 2 & 180\\
    \texttt{gen}  & 30  & 100                       & $2, 3, 4$ & 4 & 100\\
    \bottomrule
\end{tabular}
\caption{Overview of the instance sets. The generated set uses ten random job sets, each scheduled with three different machine counts.}
\label{tab:instance-sets}
\end{table}

All computational experiments were performed on an XPS 8930 workstation equipped with an Intel Core i7-8700 processor (6 physical cores and 12 logical processors, with a maximum frequency of $4.60$ GHz) and $16$ GB of RAM.
The machine ran 64-bit Ubuntu 24.04 with Linux kernel \texttt{6.17}.

\subsection{Note on computational experiments}
\label{ssec:note}
Both the instances taken from \citet{Schlembach26} and the newly generated instances represent daily sterilization operations, with all jobs intended to be completed within a single day (1440 minutes), with a time discretization of 10 minutes.

We therefore conducted all computational experiments with two alternative time horizons: the computed bound $T$ from~\eqref{eq:tmax}, and the fixed horizon $T = 144$ (in 10-minute units).
The latter reflects the daily completion requirement of the application.
For larger instances, the one-day horizon is reportedly tighter that the computed bound and yelds a smaller time-expanded network.
Testing both horizons allows us to assess the computational effect of this choice.

\subsection{Literature instances}
\label{ssec:orig-real}
The two literature sets originate from the central sterile supply department of the University Hospital of Augsburg studied by \citet{Schlembach26}.
Their data comprise ten representative working days, with observed job release dates, sizes and family assignments.
All instances have two incompatible families, with processing times $(p_1, p_2) = (6, 7)$ in 10-minute units, and capacity $q = 180$.
The latter is an integer encoding of the nine sterilization units available on each machine.
Job sizes belong to $\{1, 2, 20\}$, so that even the largest job occupies only one ninth of the machine capacity.

The \texttt{orig} set contains $150$ instances.
Each of the ten daily arrival sequences is truncated after $n$ jobs, for $n \in \{10,25,50,75,100\}$, and scheduled with $m \in \{2,3,4\}$ machines.
This gives ten instances for each of the $5 \times 3$ parameter combinations.
The smaller instances are thus prefixes of the same observed arrival sequences, and changing $m$ leaves
all job data unchanged.
This construction preserves the empirical arrival patterns while allowing the effects of problem size and machine availability to be examined systematically.

The \texttt{real} set contains the ten complete daily instances, with $m = 4$ and job counts $181$, $186$, $193$, $198$, $199$, $202$, $206$, $208$, $210$ and $214$.
Release dates range from $1$ to $83$ in the data files.
These instances retain the full workload and the machine count of the hospital setting, and complement the smaller prefix instances with tests at the scale of an actual working day.

\subsection{Generation of ``harder'' instances}
\label{ssec:gen}
The generated set complements the application-derived instances with greater variation in job sizes and family processing times, more incompatible families, and simultaneous waves of arrivals.
Its purpose is to stress the interaction between capacity-constrained batch formation and scheduling, beyond the structure of the hospital data.
We fix $n = 100$, $p = 4$ and $q = 100$, and consider $m \in \{2,3,4\}$, retaining the largest job count and the same machine counts used in \texttt{orig}.
The generator uses a reference day length of $D=144$ (in 10-minute units); this parameter scales the release-date window and is distinct from the optimization horizon $T$.\medskip

\noindent{\bf Families and processing times.}\enskip
The joint consideration of incompatible job families and unequal ready times is motivated by the semiconductor batch-scheduling setting studied by \citet{Monch05}.
One job is first assigned to each family.
Each of the remaining $n - p$ jobs is then assigned independently and uniformly to one of the $p$ families, after which the family labels are randomly shuffled.
This ensures that every family is nonempty without fixing its cardinality; each family has $n/p = 25$ jobs in expectation.
Using four families exercises the choice between incompatible batches without privileging a particular family.
Independently for each family, $p_f$ is drawn from the discrete uniform distribution over $\{2,3,\ldots,9\}$, measured in 10-minute units.
This range introduces both shorter and longer processing cycles than the fixed $60$- and $70$-minute cycles of the literature instances.
Distinct families may have the same processing time and remain incompatible.\medskip

\noindent{\bf Job sizes.}\enskip
Sampling job sizes from a discrete uniform distribution follows the experimental approach of \citet{Ozturk12}; the range is chosen here for the present stress tests.
Each size $s_j$ is drawn uniformly from $\{1,2,\ldots,\lfloor q/2 \rfloor\} = \{1,2,\ldots,50\}$, independently of its family.
The capacity $q = 100$ gives a direct percentage scale: jobs occupy between $1\%$ and $50\%$ of a machine, with an expected size of $25.5\%$ of its capacity.
This contrasts with the three size values and the maximum ratio $20/180$ in the literature sets, and makes the selection of jobs within a batch more sensitive to their sizes.
The upper limit $q/2$ ensures that any two jobs from the same family can be batched together, while the feasibility of larger batches depends on the selected combination of sizes.
The design thus retains meaningful batching choices while making capacity constraints more restrictive.\medskip

\noindent{\bf Release dates.}\enskip
The tested instances use the clustered-arrival scheme, motivated by the regular collection of medical device sets considered by \citet{Ozturk12}.
The following choice of wave count and spacing is specific to our generator.
Given a spread parameter $\gamma = 0.5$, let
\begin{equation*}
    R = \max \left\{1, \lfloor \gamma D \rfloor\right\} = 72, \qquad K = \max \left\{2, \lfloor \sqrt n \rfloor\right\} = 10.
\end{equation*}
The $K$ possible release dates are obtained by discretizing equally spaced points over $[0,R]$:
\begin{equation*}
    r_k = \left\lfloor \frac{k R}{K-1} \right\rfloor, \qquad k = 0, \ldots, K-1.
\end{equation*}
Each job independently selects one of these dates with probability $1/K$.
For the chosen parameters, arrivals therefore occur at $0, 8, 16, \ldots, 72$, with ten jobs per wave in expectation rather than a fixed wave size.
The square-root choice of $K$ provides multiple arrival waves while retaining groups of jobs that compete for the machines at the same time.
Such grouped arrivals reflect the release of instruments following surgical activity and create a tradeoff between processing available jobs immediately and waiting to form fuller batches.
Setting $\gamma = 0.5$ concentrates arrivals in the first half of the reference day and leaves the remaining half for completing the workload.
The generator controls arrivals but does not itself impose or certify completion by $D$.\medskip

\noindent{\bf Replications.}\enskip
We use ten seeds and each resulting job set is used with all three machine counts, producing a design of $30$ instances.
For a given seed, the family processing times and all job data are identical across the three values of $m$; the comparison therefore isolates the effect of machine availability.

\section{Computational results}
\label{sec:results}
We compare five implementations on the 190 instances described in Section~\ref{sec:experiments}, using both the one-day horizon $T = 144$ and the computed horizon from~\eqref{eq:tmax}, denoted here by $T_{\mathrm{UB}}$.
The comparisons address exact solution of the application-derived instances, the effect of dynamic programming on pricing, and the cost and quality of the two Lagrangian approaches.

\subsection{Compared methods and evaluation measures}
\label{ssec:measures}
The implementation labels used in the tables are defined as follows.
\begin{description}
    \item[\texttt{FLOW}] The monolithic formulation \eqref{tif:obj}--\eqref{tif:vars}, implemented in C++ with CPLEX Concert.
    \item[\texttt{CG-MIP}] The exponential-size formulation \eqref{cg:obj}--\eqref{cg:vars}, with the restricted master handled by the CPLEX C Callable Library.
    Each pricing MILP is solved exactly for a fixed family and starting time, following the decomposition used in the best exact-pricing variant of \citet{Schlembach26}.
    \item[\texttt{CG}] The same column-generation implementation, with exact pricing by the dynamic program of Section~\ref{ssec:lpcg}.
    The other algorithmic choices are the same as in \texttt{CG-MIP}.
    \item[\texttt{LAGR}] Column-generation built upon the minimum cost flow relaxation \eqref{lagr:obj}--\eqref{lagr:vars} with subgradient method of Section~\ref{ssec:cgmcf}, using the CPLEX C Callable Library and network simplex for the restricted minimum cost flow problems.
    \item[\texttt{VOL}] Column-generation built upon the shortest-path relaxation~\eqref{lagrsp:obj} with volume algorithm of Section~\ref{ssec:lagrsp}.
\end{description}
All five implementations use the initial schedule provided by the ``batch whenever idle'' heuristic of \citet{Dobson01}.
For \texttt{FLOW}, this schedule is supplied as a MILP warm start.
Both \texttt{CG-MIP} and \texttt{CG} initialize the restricted master with its batch columns and retain them in the final column pool.
\texttt{LAGR} and \texttt{VOL} also use the Dobson schedule as their starting solution.

For the four column-based methods, CPLEX solves the final integer restricted master to proven optimality within the generated column pool, without a time limit or a nonzero stopping gap.
This does not establish optimality for the full scheduling problem.
Retaining the initial schedule's columns ensures that the final integer master can reproduce that schedule; in particular, the \texttt{CG-MIP} and \texttt{CG}
solutions cannot be worse than Dobson's.
In the reported experiments, all five methods return objective values no larger than the common initial heuristic value under either horizon.
We include Dobson's solution as a baseline in the main performance tables; its standalone runtime was not recorded separately.

For instance $i$, let $z_i^A$ be the feasible objective returned by method $A$, and let $z_i^{\mathrm{ref}}$ be the smallest feasible value obtained across the five methods and both horizons.
This reference is a known optimum for all 160 literature instances and for 22 generated instances.
For the remaining eight generated instances it is a best-known value.
We measure solution quality by the relative deviation
\begin{equation}
    e_i^A = 100 \cdot \frac{z_i^A - z_i^{\mathrm{ref}}}{z_i^{\mathrm{ref}}}.
    \label{eq:results-deviation}
\end{equation}
For instances without a known optimum, this quantity measures deviation from the best observed solution, rather than from a proven optimum.
All averages are arithmetic means of the individual instance measures.

We distinguish $N_\star$, the number of solutions attaining a known optimum, from $N_c$, the number whose optimality is certified by the method's own full-problem bound.
For \texttt{FLOW}, certification follows from the branch-and-cut termination status.
For the other methods, let $lp_i^A$ be the exact master LP value for \texttt{CG-MIP} and \texttt{CG}, or the best valid full-problem Lagrangian bound for \texttt{LAGR} and \texttt{VOL}.
Since the objective is integer-valued,
\begin{equation}
    z_i^A - lp_i^A < 1
    \label{eq:results-certificate}
\end{equation}
certifies optimality.
A method can therefore attain a known optimum without certifying it.
To compare bounds independently of primal solutions, we also report
\begin{equation}
    \theta_i^A = 100 \cdot \frac{lp_i^{\texttt{CG}} - lp_i^A}{lp_i^{\texttt{CG}}},
    \label{eq:results-boundloss}
\end{equation}
using the exact column-generation LP value as a common reference.

\subsection{Results on the literature instances}
\label{ssec:literature}
Table~\ref{tab:results-literature} compares the five methods on the \texttt{orig} and \texttt{real} sets.
\texttt{FLOW} solves all 160 instances to proven optimality under both
horizons.
With the one-day horizon, its maximum runtime is 5.06 seconds on \texttt{orig} and 9.02 seconds on \texttt{real}, with means of 1.081 and 5.267 seconds, respectively.
Thus, the monolithic formulation provides an effective exact method for all literature instances in this study.
With the computed horizon, the maximum \texttt{real} runtime is 34.54 seconds.

\begin{table}[htbp]
\centering
\small
\begin{tabular}{lrrrrrrrr}
    \toprule
    & \multicolumn{4}{c}{$T = 144$} & \multicolumn{4}{c}{$T = T_{\mathrm{UB}}$}\\
    \cmidrule(lr){2-5} \cmidrule(lr){6-9}
    Method & Time & $\bar e$ & $N_\star$ & $N_c$ & Time & $\bar e$ & $N_\star$ & $N_c$\\
    \midrule
    \multicolumn{9}{@{}l}{\texttt{orig}: 150 instances, all with known optima}\\[2pt]
    %Dobson & -- & 10.125 & 15 & -- & -- & 10.125 & 15 & -- \\
    \texttt{FLOW} & 1.081 & 0.000 & 150 & 150 & 1.018 & 0.000 & 150 & 150 \\
    \texttt{CG-MIP} & 6.255 & 1.092 & 107 & 101 & 6.550 & 1.256 & 108 & 100 \\
    \texttt{CG} & 1.150 & 1.126 & 109 & 102 & 1.220 & 1.182 & 107 & 100 \\
    \texttt{LAGR} & 0.101 & 0.647 & 122 & 59 & 0.102 & 0.691 & 126 & 63 \\
    \texttt{VOL} & 0.101 & 0.016 & 145 & 130 & 0.101 & 0.016 & 145 & 130 \\
    \midrule
    \multicolumn{9}{@{}l}{\texttt{real}: 10 instances, all with known optima}\\[2pt]
    %Dobson & -- & 7.652 & 0 & -- & -- & 7.652 & 0 & -- \\
    \texttt{FLOW} & 5.267 & 0.000 & 10 & 10 & 17.632 & 0.000 & 10 & 10 \\
    \texttt{CG-MIP} & 21.514 & 0.111 & 6 & 5 & 39.209 & 0.076 & 5 & 4 \\
    \texttt{CG} & 10.230 & 0.159 & 3 & 3 & 11.539 & 0.084 & 3 & 3 \\
    \texttt{LAGR} & 0.752 & 0.000 & 10 & 0 & 0.758 & 0.007 & 9 & 0 \\
    \texttt{VOL} & 0.223 & 0.000 & 10 & 8 & 0.223 & 0.000 & 10 & 8 \\
    \bottomrule
\end{tabular}
\caption{Performance on the literature instances. Times are in seconds and $\bar e$ is the mean percentage deviation from the known optimum. $N_\star$ counts attained optima and $N_c$ counts certificates obtained by the method itself.}
\label{tab:results-literature}
\end{table}

The two exact-pricing CG variants have mean deviations from the optimum of 1.092--1.256\% on \texttt{orig} and below 0.160\% on \texttt{real}.
Depending on the variant and horizon, they certify 100--102 of the 150 \texttt{orig} instances and 3--5 of the ten \texttt{real} instances.

The Lagrangian methods substantially reduce the time needed to obtain
good feasible schedules.
In particular, \texttt{LAGR} with $T = 144$ attains all ten real-instance optima in 0.752 seconds on average, despite certifying none of them with its own bounds.
With $T = T_{\mathrm{UB}}$, it attains nine of these optima, and the remaining solution is only one objective unit above the optimum.
These outcomes illustrate why optimality attainment and certification should be reported separately: an approximate dual solution can support effective column generation even when it does not close the primal--dual gap.

\texttt{VOL} attains all ten real-instance optima under both horizons in
0.223 seconds on average, certifying eight of them.
On \texttt{orig}, it attains 145 of the 150 optima with a mean deviation of 0.016\% and certifies 130 solutions under either horizon.
\texttt{LAGR} is similarly fast on average on this set, but its mean deviation is larger.
Tables~\ref{tab:orig-details-m2}--\ref{tab:orig-details-m4} in Appendix~\ref{app:orig-results} provide the detailed breakdown of \texttt{orig} by job and machine counts.
Table~\ref{tab:real-details} in the same appendix reports averages and optimality counts over all ten \texttt{real} instances.

\subsection{Results on the generated instances}
\label{ssec:generated}
Table~\ref{tab:results-generated} separates the generated instances by machine count.
Their increased difficulty is apparent in \texttt{FLOW}: with $T = 144$, only 22 of the 30 runs certify optimality within the 600-second limit.
The counts are 3, 9 and 10 for $m = 2, 3, 4$, respectively.
With the computed horizon, these counts fall to 2, 7 and 9, for a total of 18 certificates.
The remaining runs provide feasible incumbents and valid branch-and-cut bounds.

\begin{table}[htbp]
\centering
\small
\begin{tabular}{lrrrrrrrr}
    \toprule
    & \multicolumn{4}{c}{$T = 144$} & \multicolumn{4}{c}{$T = T_{\mathrm{UB}}$}\\
    \cmidrule(lr){2-5} \cmidrule(lr){6-9}
    Method & Time & $\bar e$ & $N_\star$ & $N_c$ & Time & $\bar e$ & $N_\star$ & $N_c$\\
    \midrule
    \multicolumn{9}{@{}l}{$m=2$: 10 instances, 3 with known optima}\\[2pt]
    %Dobson & -- & 32.624 & 0 & -- & -- & 32.624 & 0 & -- \\
    \texttt{FLOW} & 459.565 & 0.112 & 3 & 3 & 501.314 & 0.404 & 3 & 2 \\
    \texttt{CG-MIP} & 14.004 & 3.832 & 0 & 0 & 46.954 & 4.242 & 0 & 0 \\
    \texttt{CG} & 18.237 & 4.178 & 0 & 0 & 16.253 & 4.273 & 0 & 0 \\
    \texttt{LAGR} & 0.724 & 6.722 & 0 & 0 & 1.110 & 8.694 & 0 & 0 \\
    \texttt{VOL} & 4.593 & 1.595 & 1 & 0 & 2.873 & 1.577 & 1 & 0 \\
    \midrule
    \multicolumn{9}{@{}l}{$m=3$: 10 instances, 9 with known optima}\\[2pt]
    %Dobson & -- & 24.009 & 0 & -- & -- & 24.009 & 0 & -- \\
    \texttt{FLOW} & 179.204 & 0.000 & 9 & 9 & 277.507 & 0.024 & 9 & 7 \\
    \texttt{CG-MIP} & 3.492 & 2.480 & 0 & 0 & 9.266 & 3.012 & 0 & 0 \\
    \texttt{CG} & 0.229 & 2.188 & 1 & 0 & 0.353 & 1.961 & 0 & 0 \\
    \texttt{LAGR} & 0.221 & 1.557 & 1 & 0 & 0.274 & 1.637 & 0 & 0 \\
    \texttt{VOL} & 0.127 & 0.504 & 3 & 0 & 0.125 & 0.504 & 3 & 0 \\
    \midrule
    \multicolumn{9}{@{}l}{$m=4$: 10 instances, 10 with known optima}\\[2pt]
    %Dobson & -- & 10.791 & 0 & -- & -- & 10.791 & 0 & -- \\
    \texttt{FLOW} & 67.259 & 0.000 & 10 & 10 & 112.020 & 0.000 & 10 & 9 \\
    \texttt{CG-MIP} & 2.781 & 1.816 & 0 & 0 & 3.530 & 1.675 & 0 & 0 \\
    \texttt{CG} & 0.078 & 1.690 & 0 & 0 & 0.098 & 1.550 & 0 & 0 \\
    \texttt{LAGR} & 0.157 & 1.482 & 1 & 0 & 0.161 & 1.296 & 0 & 0 \\
    \texttt{VOL} & 0.063 & 0.586 & 2 & 1 & 0.064 & 0.586 & 2 & 1 \\
    \bottomrule
\end{tabular}
\caption{Performance on the generated instances. Times are in seconds. $\bar e$ is the mean percentage deviation from the common best-known feasible value; this is an optimum only where optimality has been established. $N_\star$ counts matches to known optima and $N_c$ counts certificates obtained by the method itself. FLOW times include runs stopped at the 600-second limit.}
\label{tab:results-generated}
\end{table}

\texttt{FLOW} usually obtains the strongest primal result, matching the common best-known value on 29 generated instances with $T = 144$ and on 25 with $T = T_{\mathrm{UB}}$.
Its computational cost is much greater than that of the column-based heuristics.
The shortage of machines is particularly significant: the mean one-day-horizon runtime is 459.565
seconds for $m = 2$, compared with 67.259 seconds for $m = 4$.

Across the 30 generated instances, the mean deviations are 2.709\% and 2.976\% for \texttt{CG-MIP}, and 2.685\% and 2.595\% for \texttt{CG}, under the one-day and computed horizons, respectively.
The mean total times are 6.759 and 19.916 seconds for \texttt{CG-MIP}, and 6.181 and 5.568 seconds for \texttt{CG}.
Neither CG variant certifies optimality on this set; \texttt{CG} attains one known optimum with $T = 144$.

Across all 30 generated instances, \texttt{LAGR} has the smallest mean runtime: 0.368 seconds with $T = 144$ and 0.515 seconds with the computed horizon.
Its corresponding mean deviations from the common reference are 3.253\% and 3.875\%.
\texttt{VOL} obtains smaller mean deviations, 0.895\% and 0.889\%, at mean times of 1.595 and 1.021 seconds.
For three and four machines, \texttt{VOL} is also faster than the other methods on average.
For two machines, the integer optimization phase can be more expensive, as discussed below.

As an example, in an instance \texttt{VOL} with the computed horizon obtains an objective of 1692 in 11.957 seconds; the two \texttt{FLOW} runs return 1711 and 1713 after reaching the time limit.
This improves on both \texttt{FLOW} incumbents without certifying optimality.
On the generated set, \texttt{VOL} matches six known optima under either horizon, while its own bound certifies one.

\subsection{Pricing, bound computation and integer optimization}
\label{ssec:computation}
Table~\ref{tab:results-phases} separates the time spent in column generation or Lagrangian optimization (``Relax'') from the final integer-master time (``Int'').
The exact LP values returned by \texttt{CG-MIP} and \texttt{CG} coincide on all 190 instances under both horizons.
With common Dobson initialization, the comparison isolates the benefit of dynamic programming over MILP pricing.
On the generated set, the ratios of mean \texttt{CG-MIP} to mean \texttt{CG} relaxation time are 19.4 for $T = 144$ and 42.0 for the computed horizon.
These ratios concern the complete column-generation phase, including restricted-master optimization, rather than pricing time alone.
Despite their common initial solution and identical LP values, the two methods can collect different column pools and hence return different integer solutions.

\begin{table}[htbp]
\centering
\footnotesize
\begin{tabular}{llrrrrrrrr}
    \toprule
    & & \multicolumn{4}{c}{$T = 144$} & \multicolumn{4}{c}{$T = T_{\mathrm{UB}}$}\\
    \cmidrule(lr){3-6} \cmidrule(lr){7-10}
    Set & Method & Relax & Int & \#Cols & $\bar\theta$ & Relax & Int & \#Cols & $\bar\theta$\\
    \midrule
    \texttt{orig} & \texttt{CG-MIP} & 6.114 & 0.141 & 3515.2 & 0.000 & 6.393 & 0.158 & 3830.9 & 0.000 \\
     & \texttt{CG} & 1.004 & 0.147 & 3521.6 & 0.000 & 1.067 & 0.153 & 3774.2 & 0.000 \\
     & \texttt{LAGR} & 0.095 & 0.007 & 112.8 & 0.568 & 0.096 & 0.007 & 109.8 & 0.539 \\
     & \texttt{VOL} & 0.090 & 0.011 & 169.5 & 0.022 & 0.089 & 0.011 & 168.3 & 0.023 \\
    \midrule
    \texttt{real} & \texttt{CG-MIP} & 21.170 & 0.343 & 12211.7 & 0.000 & 38.956 & 0.254 & 18339.6 & 0.000 \\
     & \texttt{CG} & 9.917 & 0.313 & 12049.7 & 0.000 & 11.207 & 0.332 & 18085.2 & 0.000 \\
     & \texttt{LAGR} & 0.738 & 0.014 & 330.3 & 1.237 & 0.742 & 0.015 & 336.4 & 1.196 \\
     & \texttt{VOL} & 0.208 & 0.015 & 339.5 & 0.023 & 0.208 & 0.015 & 339.5 & 0.023 \\
    \midrule
    \texttt{gen} & \texttt{CG-MIP} & 3.899 & 2.860 & 2612.4 & 0.000 & 12.415 & 7.501 & 3568.0 & 0.000 \\
     & \texttt{CG} & 0.201 & 5.981 & 2558.2 & 0.000 & 0.296 & 5.273 & 3871.5 & 0.000 \\
     & \texttt{LAGR} & 0.205 & 0.163 & 750.3 & 1.656 & 0.282 & 0.233 & 800.3 & 1.605 \\
     & \texttt{VOL} & 0.103 & 1.492 & 639.3 & 0.030 & 0.105 & 0.915 & 639.7 & 0.029 \\
    \bottomrule
\end{tabular}
\caption{Cost of the relaxation and integer phases. Times are in seconds, ``\#Cols'' is the mean reported column count, and $\bar\theta$ is the mean percentage lower-bound loss relative to the exact CG relaxation.}
\label{tab:results-phases}
\end{table}

The two Lagrangian procedures report substantially fewer columns.
For example, with $T = 144$, the mean counts on \texttt{real} are 12050 for \texttt{CG}, 330 for \texttt{LAGR} and 340 for \texttt{VOL}.
\texttt{VOL} also obtains bounds close to the exact LP values: its mean relative bound losses are 0.022\%, 0.023\% and 0.030\% on \texttt{orig}, \texttt{real} and \texttt{gen}, respectively.
The corresponding losses for \texttt{LAGR} are larger, although its generated columns still yield competitive integer solutions, especially on the real-world instances.

Fast relaxation does not always imply a similarly fast complete algorithm.
On \texttt{gen} with $T = 144$, \texttt{CG} spends 0.201 seconds on its relaxation phase and 5.981 seconds on integer optimization, on average.
\texttt{CG-MIP} spends 3.899 and 2.860 seconds, respectively.
Thus, the large relaxation-phase speedup yields only a modest reduction in mean total time under this horizon.
For \texttt{VOL}, the phase times are 0.103 and 1.492 seconds.
Its generated-instance median total runtime is only 0.096 seconds, so a small number of difficult integer masters have a substantial effect on the mean.

\subsection{Effect of the time horizon}
\label{ssec:horizon}
Table~\ref{tab:results-horizons} relates the computed horizon to the runtime differences.
The computed horizon is shorter than 144 for all 10- and 25-job instances, but exceeds 144 for all ten real-world instances and for all generated instances with two machines.
It is shorter for 113 of the 150 \texttt{orig} instances and for six of the 30 generated instances.
The common exact CG LP value is unchanged between the two horizons on every instance; changes in runtime and integer solutions therefore occur without a change in this relaxation value.

\begin{table}[htbp]
\centering
\footnotesize
\begin{tabular}{lrrr@{\,}c@{\,}lrrrrrr}
    \toprule
    & & & & & & & \multicolumn{5}{c}{Runtime ratio: UB / T}\\
    \cmidrule(lr){8-12}
    Group & $N$ & $\bar T_{\mathrm{UB}}$ & \multicolumn{3}{c}{Range} & $N_{<}$ & \texttt{FLOW} & \texttt{CG-MIP} & \texttt{CG} & \texttt{LAGR} & \texttt{VOL}\\
    \midrule
    \texttt{orig}, $n=10$ & 30 & 16.9 & 7&--&35 & 30 & 0.34 & 0.11 & 0.57 & 1.14 & 1.03 \\
    \texttt{orig}, $n=25$ & 30 & 50.7 & 19&--&118 & 30 & 0.51 & 0.33 & 0.90 & 0.97 & 0.98 \\
    \texttt{orig}, $n=50$ & 30 & 100.3 & 36&--&202 & 24 & 0.62 & 0.76 & 0.97 & 1.04 & 0.99 \\
    \texttt{orig}, $n=75$ & 30 & 144.9 & 67&--&278 & 19 & 0.97 & 1.10 & 1.00 & 0.99 & 1.00 \\
    \texttt{orig}, $n=100$ & 30 & 205.1 & 107&--&433 & 10 & 1.07 & 1.44 & 1.11 & 1.02 & 1.00 \\
    \texttt{real} & 10 & 268.6 & 223&--&303 & 0 & 3.35 & 1.82 & 1.13 & 1.01 & 1.00 \\
    \texttt{gen}, $m=2$ & 10 & 845.4 & 276&--&1677 & 0 & 1.09 & 3.35 & 0.89 & 1.53 & 0.63 \\
    \texttt{gen}, $m=3$ & 10 & 380.9 & 138&--&754 & 1 & 1.55 & 2.65 & 1.54 & 1.24 & 0.98 \\
    \texttt{gen}, $m=4$ & 10 & 187.3 & 107&--&401 & 5 & 1.67 & 1.27 & 1.27 & 1.02 & 1.00 \\
    \bottomrule
\end{tabular}
\caption{Computed horizons and runtime ratios. Average horizons $\bar T_{\mathrm{UB}}$ and ranges are in minutes; $N_{<}$ counts computed horizons shorter than one day ($T = 144$). Each ratio is the sum of UB-run times divided by the sum of T-run times in the same group. FLOW ratios include time-limited runs.}
\label{tab:results-horizons}
\end{table}

The tighter computed horizons reduce the mean times of \texttt{FLOW} and the two exact-pricing CG methods on the 10- and 25-job instances.
The one-day limit can be computationally useful for larger workloads.
For example, on \texttt{real}, using the computed horizon increases mean runtime by factors of 3.35 for \texttt{FLOW}, 1.82 for \texttt{CG-MIP}, and 1.13 for \texttt{CG}.
On the generated instances with two machines, the \texttt{CG-MIP} ratio is 3.35, whereas the \texttt{CG} ratio is 0.89.
These aggregate ratios show that the horizon does not affect all methods uniformly: it can alter the generated pool and the subsequent integer optimization.
For \texttt{FLOW}, ratios involving time-limited runs compare consumed runtime budgets rather than times to proven optimality.

\texttt{LAGR} and \texttt{VOL} are less sensitive to the horizon on the literature instances.
On \texttt{gen}, the lower mean time of \texttt{VOL} with the computed horizon comes from the integer phase.
A known daily completion limit is valuable when available, while \eqref{eq:tmax} supplies a safe horizon when no such limit is known.

\section{Concluding remarks}
\label{sec:conclusions}

We have studied total flow time minimization on identical parallel
batch-processing machines with release dates, unequal job sizes and
incompatible families, motivated by hospital sterilization operations.
The proposed arc-flow formulations provide a common representation for a
pseudopolynomial mixed-integer model, column generation and two Lagrangian
relaxations.

The Lagrangian approaches generate substantially smaller column pools and
obtain good feasible schedules at low computational cost. In particular,
the shortest-path relaxation combined with the volume algorithm offers a
favorable balance between solution quality, lower-bound accuracy and
computational effort.

Further research should address the recovery of integer solutions from the
generated batches. Rounding and repair procedures based on the primal
averages produced by the volume algorithm could reduce the cost of the
final integer optimization.

The framework also suggests extensions to related batch scheduling problems
with weighted flow-time or tardiness objectives
--- along the lines of \citet{Pessoa22} --- with heterogeneous machines,
machine availability restrictions and family-dependent setup times.
A central question is which coupling constraints can be relaxed while
retaining efficient network subproblems and effective batch generation.

\clearpage
\appendix

\numberwithin{equation}{section}
\counterwithin{table}{section}
\counterwithin{figure}{section}

\section{Detailed results for the literature benchmark}
\label{app:orig-results}
Tables~\ref{tab:orig-details-m2}--\ref{tab:orig-details-m4} report the complete breakdown of \texttt{orig} by job and machine counts.
Table~\ref{tab:real-details} summarizes all ten \texttt{real} instances with $m=4$.
In all four tables, each row reports averages over ten instances and the counts $(N_\star,N_c)$.
The deviation $\bar e$ uses the known optimum, and $(N_\star,N_c)$ distinguishes attaining that optimum from certifying it.
We additionally report the mean relative primal--dual gap $\bar g$: it is set to zero when the method certifies optimality; otherwise the reported percentage $100 (z_i^A - lp_i^A) / z_i^A$ is used.
For \texttt{FLOW}, the bound is the final branch-and-cut bound, not the value of the initial LP relaxation.

\begin{table}[htbp]
\centering
\footnotesize
\begin{tabular}{rlrrrcrrrc}
    \toprule
    & & \multicolumn{4}{c}{$T = 144$} & \multicolumn{4}{c}{$T = T_{\mathrm{UB}}$}\\
    \cmidrule(lr){3-6} \cmidrule(lr){7-10}
    $n$ & Method & Time & $\bar e$ & $\bar g$ & $(N_\star,N_c)$ & Time & $\bar e$ & $\bar g$ & $(N_\star,N_c)$\\
    \midrule
    10 & \texttt{FLOW} & 0.038 & 0.000 & 0.000 & $(10,10)$ & 0.016 & 0.000 & 0.000 & $(10,10)$ \\
     & \texttt{CG-MIP} & 0.924 & 0.000 & 0.000 & $(10,10)$ & 0.156 & 0.000 & 0.000 & $(10,10)$ \\
     & \texttt{CG} & 0.009 & 0.000 & 0.000 & $(10,10)$ & 0.007 & 0.000 & 0.000 & $(10,10)$ \\
     & \texttt{LAGR} & 0.011 & 0.000 & 0.000 & $(10,10)$ & 0.012 & 0.000 & 0.000 & $(10,10)$ \\
     & \texttt{VOL} & 0.014 & 0.000 & 0.000 & $(10,10)$ & 0.013 & 0.000 & 0.000 & $(10,10)$ \\
    \midrule
    25 & \texttt{FLOW} & 0.170 & 0.000 & 0.000 & $(10,10)$ & 0.132 & 0.000 & 0.000 & $(10,10)$ \\
     & \texttt{CG-MIP} & 4.647 & 2.002 & 2.082 & $(8,6)$ & 2.231 & 2.418 & 2.437 & $(8,6)$ \\
     & \texttt{CG} & 0.064 & 1.468 & 1.639 & $(8,6)$ & 0.063 & 1.311 & 1.548 & $(8,6)$ \\
     & \texttt{LAGR} & 0.045 & 1.703 & 1.863 & $(9,7)$ & 0.043 & 1.703 & 1.898 & $(9,7)$ \\
     & \texttt{VOL} & 0.068 & 0.000 & 0.328 & $(10,7)$ & 0.068 & 0.000 & 0.328 & $(10,7)$ \\
    \midrule
    50 & \texttt{FLOW} & 0.604 & 0.000 & 0.000 & $(10,10)$ & 0.631 & 0.000 & 0.000 & $(10,10)$ \\
     & \texttt{CG-MIP} & 8.427 & 4.145 & 3.771 & $(5,5)$ & 8.749 & 3.860 & 3.544 & $(6,6)$ \\
     & \texttt{CG} & 1.161 & 4.629 & 4.148 & $(6,6)$ & 1.331 & 5.078 & 4.645 & $(5,5)$ \\
     & \texttt{LAGR} & 0.115 & 3.555 & 3.427 & $(8,4)$ & 0.109 & 3.920 & 3.693 & $(7,7)$ \\
     & \texttt{VOL} & 0.161 & 0.055 & 0.277 & $(9,8)$ & 0.160 & 0.055 & 0.277 & $(9,8)$ \\
    \midrule
    75 & \texttt{FLOW} & 1.145 & 0.000 & 0.000 & $(10,10)$ & 2.048 & 0.000 & 0.000 & $(10,10)$ \\
     & \texttt{CG-MIP} & 10.017 & 3.801 & 3.452 & $(4,4)$ & 15.771 & 5.286 & 4.768 & $(3,3)$ \\
     & \texttt{CG} & 2.498 & 3.959 & 3.599 & $(4,4)$ & 2.863 & 3.650 & 3.367 & $(4,4)$ \\
     & \texttt{LAGR} & 0.151 & 2.180 & 2.422 & $(4,2)$ & 0.166 & 1.834 & 2.134 & $(6,2)$ \\
     & \texttt{VOL} & 0.238 & 0.000 & 0.054 & $(10,8)$ & 0.239 & 0.000 & 0.054 & $(10,8)$ \\
    \midrule
    100 & \texttt{FLOW} & 2.675 & 0.000 & 0.000 & $(10,10)$ & 3.817 & 0.000 & 0.000 & $(10,10)$ \\
     & \texttt{CG-MIP} & 14.888 & 4.849 & 4.587 & $(2,2)$ & 29.087 & 5.766 & 5.468 & $(1,1)$ \\
     & \texttt{CG} & 4.738 & 5.268 & 5.027 & $(1,1)$ & 5.730 & 6.185 & 5.752 & $(0,0)$ \\
     & \texttt{LAGR} & 0.253 & 1.449 & 2.027 & $(3,1)$ & 0.261 & 1.931 & 2.569 & $(2,0)$ \\
     & \texttt{VOL} & 0.373 & 0.157 & 0.360 & $(7,5)$ & 0.373 & 0.157 & 0.360 & $(7,5)$ \\
    \bottomrule
\end{tabular}
\caption{Detailed results for \texttt{orig} with $m=2$; each row summarizes ten instances. Times are in seconds; $\bar e$ is mean deviation from the optimum, $\bar g$ is the mean primal--dual gap (zero when certified), and the last entry is the pair $(N_\star,N_c)$.}
\label{tab:orig-details-m2}
\end{table}

\begin{table}[htbp]
\centering
\footnotesize
\begin{tabular}{rlrrrcrrrc}
    \toprule
    & & \multicolumn{4}{c}{$T = 144$} & \multicolumn{4}{c}{$T = T_{\mathrm{UB}}$}\\
    \cmidrule(lr){3-6} \cmidrule(lr){7-10}
    $n$ & Method & Time & $\bar e$ & $\bar g$ & $(N_\star,N_c)$ & Time & $\bar e$ & $\bar g$ & $(N_\star,N_c)$\\
    \midrule
    10 & \texttt{FLOW} & 0.036 & 0.000 & 0.000 & $(10,10)$ & 0.014 & 0.000 & 0.000 & $(10,10)$ \\
     & \texttt{CG-MIP} & 0.674 & 0.000 & 0.000 & $(10,10)$ & 0.044 & 0.000 & 0.000 & $(10,10)$ \\
     & \texttt{CG} & 0.011 & 0.000 & 0.000 & $(10,10)$ & 0.006 & 0.000 & 0.000 & $(10,10)$ \\
     & \texttt{LAGR} & 0.012 & 0.000 & 0.000 & $(10,10)$ & 0.017 & 0.000 & 0.000 & $(10,10)$ \\
     & \texttt{VOL} & 0.010 & 0.000 & 0.000 & $(10,10)$ & 0.011 & 0.000 & 0.000 & $(10,10)$ \\
    \midrule
    25 & \texttt{FLOW} & 0.219 & 0.000 & 0.000 & $(10,10)$ & 0.118 & 0.000 & 0.000 & $(10,10)$ \\
     & \texttt{CG-MIP} & 3.669 & 0.119 & 0.162 & $(9,9)$ & 0.999 & 0.119 & 0.162 & $(9,9)$ \\
     & \texttt{CG} & 0.036 & 0.119 & 0.162 & $(9,9)$ & 0.030 & 0.119 & 0.162 & $(9,9)$ \\
     & \texttt{LAGR} & 0.030 & 0.119 & 0.438 & $(9,7)$ & 0.031 & 0.119 & 0.443 & $(9,7)$ \\
     & \texttt{VOL} & 0.037 & 0.000 & 0.000 & $(10,10)$ & 0.036 & 0.000 & 0.000 & $(10,10)$ \\
    \midrule
    50 & \texttt{FLOW} & 0.767 & 0.000 & 0.000 & $(10,10)$ & 0.423 & 0.000 & 0.000 & $(10,10)$ \\
     & \texttt{CG-MIP} & 5.762 & 0.181 & 0.261 & $(8,7)$ & 3.583 & 0.031 & 0.138 & $(9,7)$ \\
     & \texttt{CG} & 0.440 & 0.509 & 0.565 & $(9,8)$ & 0.304 & 0.123 & 0.230 & $(8,6)$ \\
     & \texttt{LAGR} & 0.070 & 0.361 & 1.290 & $(6,1)$ & 0.081 & 0.534 & 1.380 & $(8,1)$ \\
     & \texttt{VOL} & 0.073 & 0.000 & 0.085 & $(10,8)$ & 0.072 & 0.000 & 0.085 & $(10,8)$ \\
    \midrule
    75 & \texttt{FLOW} & 1.623 & 0.000 & 0.000 & $(10,10)$ & 1.344 & 0.000 & 0.000 & $(10,10)$ \\
     & \texttt{CG-MIP} & 8.455 & 0.315 & 0.407 & $(8,6)$ & 7.155 & 0.414 & 0.499 & $(8,6)$ \\
     & \texttt{CG} & 1.483 & 0.199 & 0.312 & $(7,5)$ & 1.299 & 0.178 & 0.275 & $(8,6)$ \\
     & \texttt{LAGR} & 0.138 & 0.159 & 0.992 & $(9,1)$ & 0.133 & 0.159 & 1.042 & $(9,0)$ \\
     & \texttt{VOL} & 0.119 & 0.000 & 0.099 & $(10,7)$ & 0.119 & 0.000 & 0.099 & $(10,7)$ \\
    \midrule
    100 & \texttt{FLOW} & 2.578 & 0.000 & 0.000 & $(10,10)$ & 2.682 & 0.000 & 0.000 & $(10,10)$ \\
     & \texttt{CG-MIP} & 11.388 & 0.324 & 0.387 & $(5,4)$ & 14.628 & 0.382 & 0.452 & $(5,3)$ \\
     & \texttt{CG} & 3.410 & 0.150 & 0.203 & $(7,5)$ & 3.540 & 0.263 & 0.315 & $(6,5)$ \\
     & \texttt{LAGR} & 0.218 & 0.148 & 1.078 & $(6,0)$ & 0.210 & 0.148 & 0.988 & $(7,0)$ \\
     & \texttt{VOL} & 0.165 & 0.029 & 0.085 & $(9,7)$ & 0.167 & 0.029 & 0.085 & $(9,7)$ \\
    \bottomrule
\end{tabular}
\caption{Detailed results for \texttt{orig} with $m=3$; each row summarizes ten instances. Times are in seconds; $\bar e$ is mean deviation from the optimum, $\bar g$ is the mean primal--dual gap (zero when certified), and the last entry is the pair $(N_\star,N_c)$.}
\label{tab:orig-details-m3}
\end{table}

\begin{table}[htbp]
\centering
\footnotesize
\begin{tabular}{rlrrrcrrrc}
    \toprule
    & & \multicolumn{4}{c}{$T = 144$} & \multicolumn{4}{c}{$T = T_{\mathrm{UB}}$}\\
    \cmidrule(lr){3-6} \cmidrule(lr){7-10}
    $n$ & Method & Time & $\bar e$ & $\bar g$ & $(N_\star,N_c)$ & Time & $\bar e$ & $\bar g$ & $(N_\star,N_c)$\\
    \midrule
    10 & \texttt{FLOW} & 0.043 & 0.000 & 0.000 & $(10,10)$ & 0.010 & 0.000 & 0.000 & $(10,10)$ \\
     & \texttt{CG-MIP} & 0.503 & 0.000 & 0.000 & $(10,10)$ & 0.024 & 0.000 & 0.000 & $(10,10)$ \\
     & \texttt{CG} & 0.016 & 0.000 & 0.000 & $(10,10)$ & 0.007 & 0.000 & 0.000 & $(10,10)$ \\
     & \texttt{LAGR} & 0.012 & 0.000 & 0.000 & $(10,10)$ & 0.012 & 0.000 & 0.000 & $(10,10)$ \\
     & \texttt{VOL} & 0.008 & 0.000 & 0.000 & $(10,10)$ & 0.009 & 0.000 & 0.000 & $(10,10)$ \\
    \midrule
    25 & \texttt{FLOW} & 0.296 & 0.000 & 0.000 & $(10,10)$ & 0.100 & 0.000 & 0.000 & $(10,10)$ \\
     & \texttt{CG-MIP} & 2.880 & 0.183 & 0.220 & $(9,9)$ & 0.519 & 0.122 & 0.161 & $(9,9)$ \\
     & \texttt{CG} & 0.027 & 0.000 & 0.000 & $(10,10)$ & 0.020 & 0.122 & 0.161 & $(9,9)$ \\
     & \texttt{LAGR} & 0.040 & 0.000 & 0.648 & $(10,6)$ & 0.037 & 0.000 & 0.318 & $(10,7)$ \\
     & \texttt{VOL} & 0.024 & 0.000 & 0.000 & $(10,10)$ & 0.023 & 0.000 & 0.000 & $(10,10)$ \\
    \midrule
    50 & \texttt{FLOW} & 0.995 & 0.000 & 0.000 & $(10,10)$ & 0.401 & 0.000 & 0.000 & $(10,10)$ \\
     & \texttt{CG-MIP} & 5.001 & 0.064 & 0.085 & $(8,8)$ & 2.320 & 0.064 & 0.085 & $(8,8)$ \\
     & \texttt{CG} & 0.268 & 0.064 & 0.085 & $(8,8)$ & 0.185 & 0.318 & 0.331 & $(8,8)$ \\
     & \texttt{LAGR} & 0.075 & 0.000 & 1.088 & $(10,0)$ & 0.079 & 0.000 & 0.859 & $(10,2)$ \\
     & \texttt{VOL} & 0.045 & 0.000 & 0.000 & $(10,10)$ & 0.045 & 0.000 & 0.000 & $(10,10)$ \\
    \midrule
    75 & \texttt{FLOW} & 2.015 & 0.000 & 0.000 & $(10,10)$ & 1.240 & 0.000 & 0.000 & $(10,10)$ \\
     & \texttt{CG-MIP} & 6.984 & 0.147 & 0.171 & $(5,5)$ & 5.000 & 0.169 & 0.192 & $(6,6)$ \\
     & \texttt{CG} & 0.972 & 0.273 & 0.295 & $(5,5)$ & 0.767 & 0.212 & 0.224 & $(6,6)$ \\
     & \texttt{LAGR} & 0.147 & 0.000 & 1.098 & $(10,0)$ & 0.132 & 0.000 & 1.243 & $(10,0)$ \\
     & \texttt{VOL} & 0.073 & 0.000 & 0.000 & $(10,10)$ & 0.073 & 0.000 & 0.000 & $(10,10)$ \\
    \midrule
    100 & \texttt{FLOW} & 3.005 & 0.000 & 0.000 & $(10,10)$ & 2.299 & 0.000 & 0.000 & $(10,10)$ \\
     & \texttt{CG-MIP} & 9.609 & 0.249 & 0.256 & $(6,6)$ & 7.986 & 0.204 & 0.212 & $(6,6)$ \\
     & \texttt{CG} & 2.125 & 0.247 & 0.255 & $(5,5)$ & 2.147 & 0.170 & 0.187 & $(6,6)$ \\
     & \texttt{LAGR} & 0.200 & 0.031 & 0.982 & $(8,0)$ & 0.213 & 0.015 & 1.156 & $(9,0)$ \\
     & \texttt{VOL} & 0.103 & 0.000 & 0.000 & $(10,10)$ & 0.103 & 0.000 & 0.000 & $(10,10)$ \\
    \bottomrule
\end{tabular}
\caption{Detailed results for \texttt{orig} with $m=4$; each row summarizes ten instances. Times are in seconds; $\bar e$ is mean deviation from the optimum, $\bar g$ is the mean primal--dual gap (zero when certified), and the last entry is the pair $(N_\star,N_c)$.}
\label{tab:orig-details-m4}
\end{table}

\begin{table}[htbp]
\centering
\footnotesize
\begin{tabular}{lrrrcrrrc}
    \toprule
    & \multicolumn{4}{c}{$T = 144$} & \multicolumn{4}{c}{$T = T_{\mathrm{UB}}$}\\
    \cmidrule(lr){2-5} \cmidrule(lr){6-9}
    Method & Time & $\bar e$ & $\bar g$ & $(N_\star,N_c)$ & Time & $\bar e$ & $\bar g$ & $(N_\star,N_c)$\\
    \midrule
    \texttt{FLOW} & 5.267 & 0.000 & 0.000 & $(10,10)$ & 17.632 & 0.000 & 0.000 & $(10,10)$ \\
    \texttt{CG-MIP} & 21.514 & 0.111 & 0.130 & $(6,5)$ & 39.209 & 0.076 & 0.095 & $(5,4)$ \\
    \texttt{CG} & 10.230 & 0.159 & 0.178 & $(3,3)$ & 11.539 & 0.084 & 0.103 & $(3,3)$ \\
    \texttt{LAGR} & 0.752 & 0.000 & 1.256 & $(10,0)$ & 0.758 & 0.007 & 1.222 & $(9,0)$ \\
    \texttt{VOL} & 0.223 & 0.000 & 0.024 & $(10,8)$ & 0.223 & 0.000 & 0.024 & $(10,8)$ \\
    \bottomrule
\end{tabular}
\caption{Detailed results for \texttt{real} with $m=4$; each row summarizes all ten instances. Times are in seconds; $\bar e$ is mean deviation from the optimum, $\bar g$ is the mean primal--dual gap (zero when certified), and the last entry is the pair $(N_\star,N_c)$.}
\label{tab:real-details}
\end{table}

\clearpage

\section*{CRediT authorship contribution statement}
\textbf{Alessandro Druetto:} Conceptualization, Methodology, Software, Validation, Writing - Original Draft, Writing - Review \& Editing.
\textbf{Andrea Grosso:} Conceptualization, Methodology, Software, Validation, Writing - Original Draft, Writing - Review \& Editing.

\section*{Funding}
This research did not receive any specific grant from funding agencies
in the public, commercial, or not-for-profit sectors.

\section*{Declaration of competing interest}
The authors do not recognize any conflicting interest.

\section*{Data availability}
The literature benchmark instances originate from \citet{Schlembach26}.
Generated instances and algorithm implementations are available upon request to the authors.

\section*{Declaration of generative AI and AI-assisted technologies in the manuscript preparation process}
During the preparation of this work, the authors used Claude (Anthropic) to help in code development and ChatGPT (OpenAI) to assist with drafting and revising the manuscript.
After using these tools, the authors reviewed and edited the content as needed and take full responsibility for the content of the published article.

\section*{Acknowledgements}
Authors are grateful to Proff.\ J.~O.~Brunner and R.~Schlembach for having provided useful 
clarifications about~\citep{Schlembach26} in the early stages of this research.

\end{document}